\documentclass[article]{elsarticle}

\usepackage{lineno,hyperref}
\usepackage{amsmath,amssymb,graphicx}
\usepackage{mathtools}
\usepackage{engord}
\usepackage{color}

\usepackage{amsthm,bm}
\usepackage{thmtools}
\usepackage{mathrsfs}
\usepackage{subfigure}

\newtheorem{proposition}{\textbf{Proposition}}

\theoremstyle{definition}
\newtheorem{definition}{\textbf{Definition}}

\declaretheoremstyle[%
  spaceabove=-6pt,%
  spacebelow=6pt,%
  headfont=\normalfont\itshape,%
  postheadspace=1em,%
  qed=\qedsymbol,%
  headpunct={}
]{mystyle}

\def\R{\mathbb{R}}
\def\Z{\mathbb{Z}}

\begin{document}

\begin{frontmatter}

\title{Compactly-Supported Smooth Interpolators for Shape Modeling with Varying Resolution}




\author[mymainaddress]{D. Schmitter, J. Fageot, A. Badoual, P. Garcia-Amorena, and M. Unser}



\address[mymainaddress]{Biomedical Imaging Group, EPFL, Switzerland}

\begin{abstract}
	In applications that involve interactive curve and surface modeling, the intuitive manipulation of shapes is crucial. For instance, user interaction is facilitated if a geometrical object can be manipulated through control points that interpolate the shape itself. Additionally, models for shape representation often need to provide local shape control and they need to be able to reproduce common shape primitives such as ellipsoids, spheres, cylinders, or tori. We present a general framework to construct families of compactly-supported interpolators that are piecewise-exponential polynomial. They can be designed to satisfy regularity constraints of any order and they enable one to build parametric deformable shape models by suitable linear combinations of interpolators.  They allow to change the resolution of shapes based on the refinability of B-splines. We illustrate their use on examples to construct shape models that involve curves and surfaces with applications to interactive modeling and character design.

\end{abstract}


\end{frontmatter}


\section{Introduction}
\noindent The interactive modeling of curves and surfaces is desirable in applications that involve the visualization of shapes. Related domains include computer graphics~\cite{Brechbuhler1995,Romani2004,Conti2015bis,Dyn1987,Cohen1980,botsch2010}, image analysis in biomedical imaging~\cite{Delgado13B,Tang2004,Dufour2011,Barbosa2012,Schmitter14}, industrial shape design~\cite{Audette2012,Garg2014,Deng2015} or the modeling of animated surfaces~\cite{DeRose1998}. Shape-modeling frameworks that allow for user interaction can usually be categorized in either \textit{discrete} or \textit{continuous-domain} models. Discrete models are typically based on interpolating polygon meshes or subdivision~\cite{Dyn2002,Dyn2007,Conti2011,Charina2014,Novara2015,Badoual2017} and they easily allow to locally refine a shape. Subdivision models are also considered as hybrids between discrete and continous-domain models because they iteratively define continuous functions in the limit. However, the limit functions do not always have a closed-form expression~\cite{Dyn2003}. Continuous-domain models allow for \textit{organic} shape modeling and consist of B\'ezier shapes or spline-based models such as NURBS~\cite{Piegl1991,Piegl2010,Rogers2001}. They allow one to control shapes locally due to their {compactly} supported basis functions. However, NURBS generally cannot be smooth and interpolating at the same time, which leads to a non-intuitive manipulation of shapes because NURBS control points do not lie on the boundary of the object. 

\subsection{Motivation and Contribution}
\begin{figure}[htb]
\begin{minipage}{\textwidth}
\begin{centering}
\includegraphics[width=11.6cm, trim=0.6cm 0.6cm 0.6cm 0.6cm, clip=true]{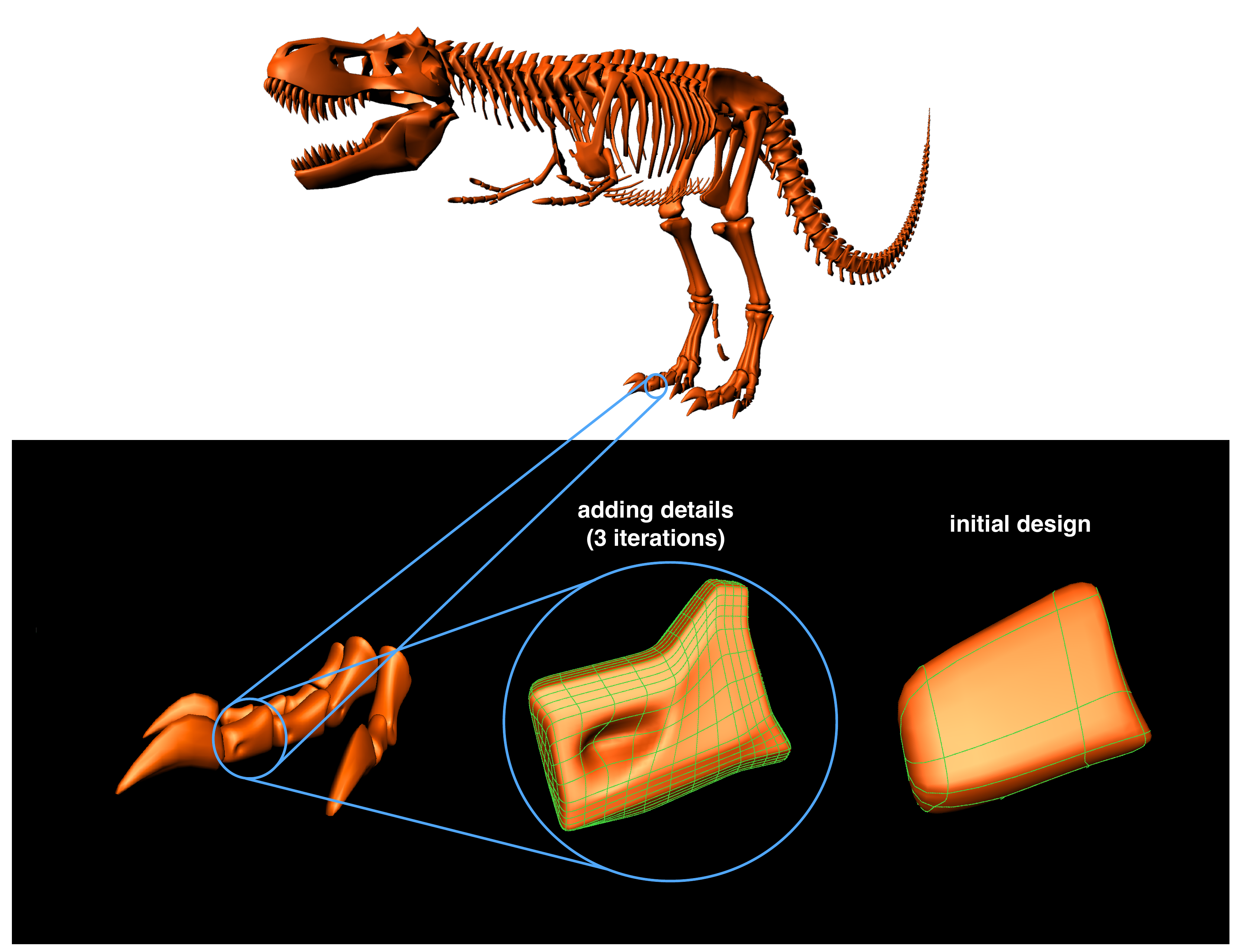}
\caption[The LOF caption]{Interactive shape modeling for character design. Remodeling of the foot of the ``T-rex'' is shown. A bone of the middle toe of the right foot is modeled; first, an initial design is achieved with few control points that interpolate the shape (bottom, right). Then, the resolution is increased by applying three refinement iterations in order to have more flexibility to add details to the bone (bottom, middle). Due to convergence of our modified refinement scheme, after three iterations it behaves interpolatory-like.}
\label{fig: fancy shapes}
\end{centering}
\footnotetext{\tiny{The ``T-rex'' has been remodeled after the character designed by Joel Anderson, source: http://joel3d.com/}}
\end{minipage}
\end{figure}

\begin{figure}[!ht]
\center
\subfigure[Torus.]{
\includegraphics[width=3.6cm, trim=1cm 1cm 1cm 1cm, clip=true]{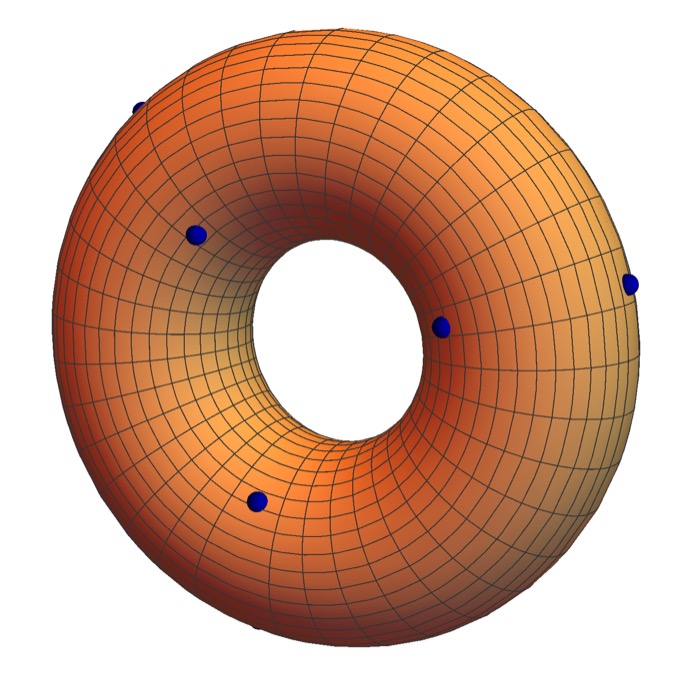}
}\hspace*{+1cm}
\subfigure[``Figure 8'' immersion.]{
\includegraphics[width=3.6cm, trim=1cm 1cm 1cm 1cm, clip=true]{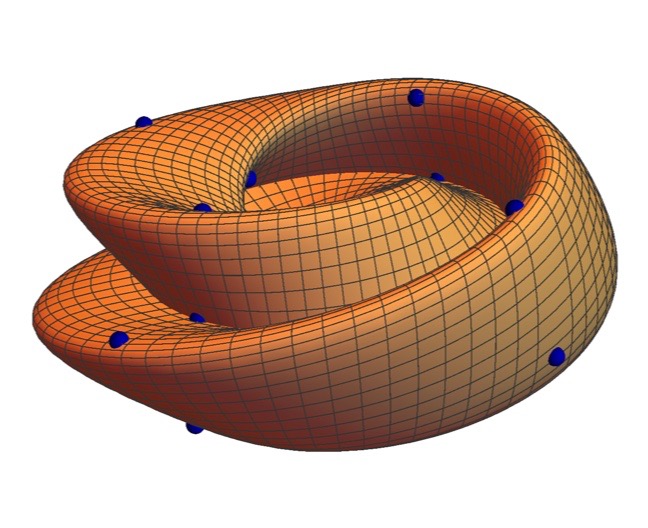}
}
\\
\subfigure[Helicoid.]{
\includegraphics[width=4.2cm, trim=1cm 2.5cm 1cm 1cm, clip=true]{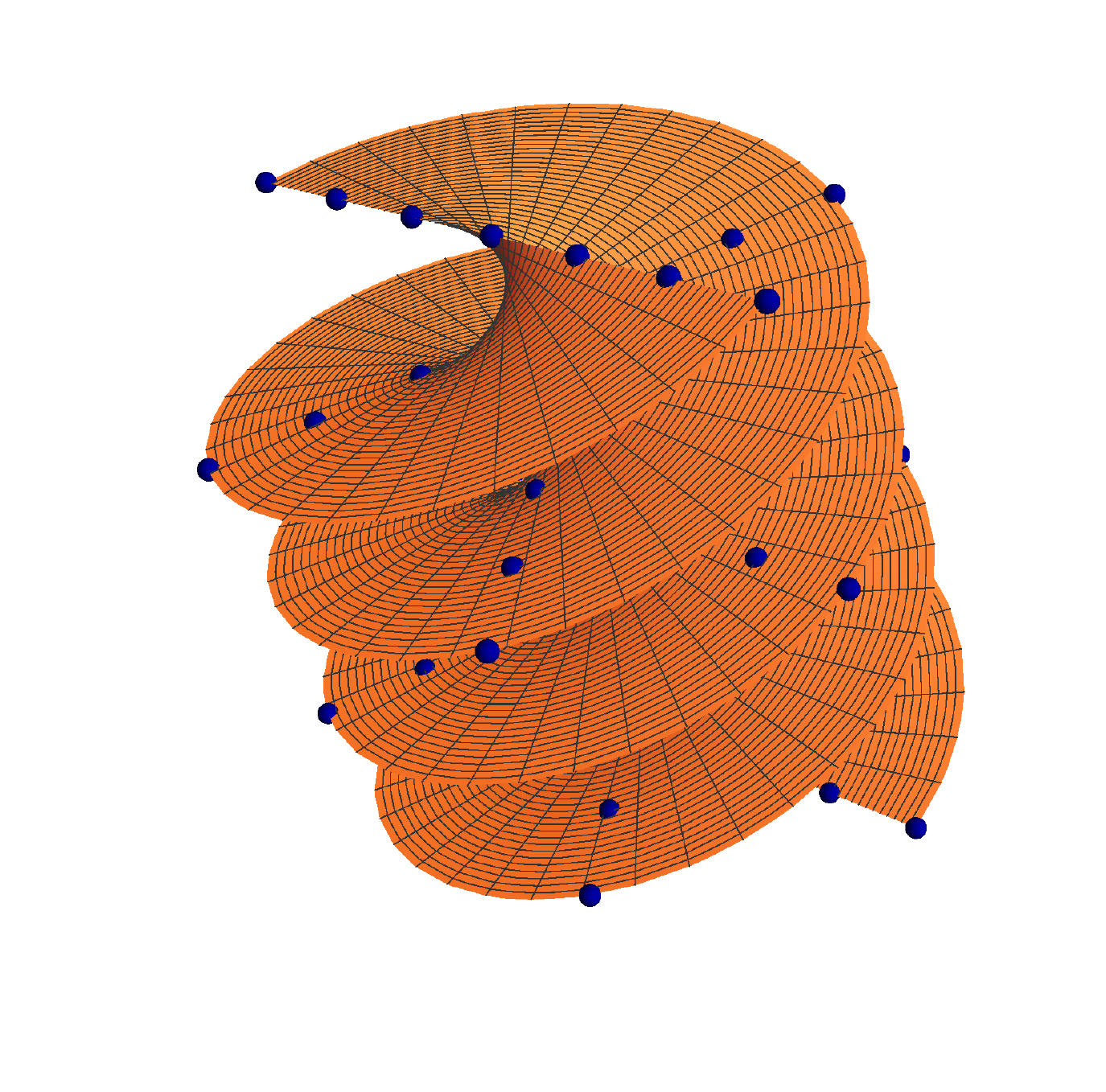}
}\hspace*{+1cm}
\subfigure[Pinched torus.]{
\includegraphics[width=3.6cm, trim=1cm 3cm 1cm 1cm, clip=true]{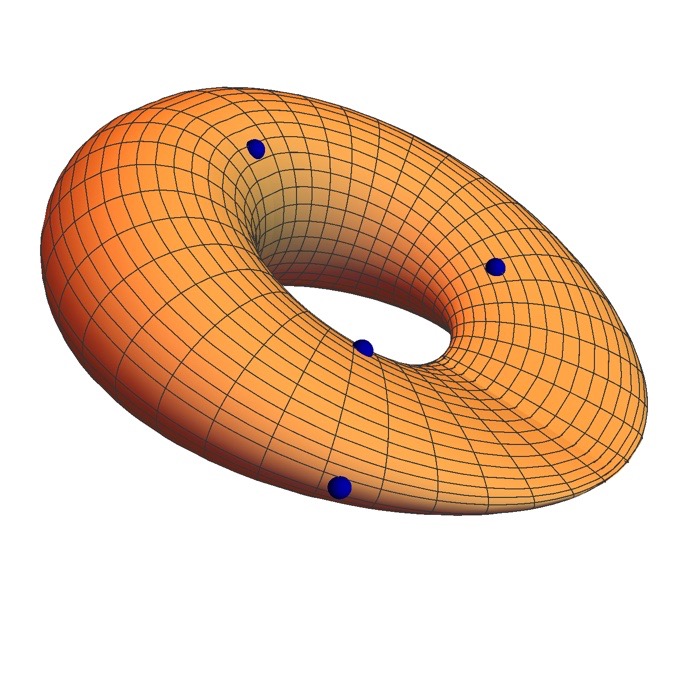}
}
\caption{Parametric surfaces constructed with the proposed family of interpolators. If the parameterization of a shape is known, we provide the formulae to construct the corresponding interpolator in order to represent the shape as detailed in Section~\ref{sec: applications}. The interpolation property ensures that the control points (blue points) interpolate the surface. This property is particularly useful in user-interactive applications, where a surface is modified by dragging control points (\textit{e.g.} as previously demonstrated in~\cite{schmitter15,schmitter16b})}
\label{fig:illustrations}
\end{figure}

\noindent Our motivation is the practical need for interpolating functions to be used in user-interactive applications\footnote{Videos that illustrate the use and advantage of our proposed framework can be found at http://bigwww.epfl.ch/demo/varying-resolution-interpolator/.} (see Figures~\ref{fig: fancy shapes} and~\ref{fig:illustrations}). In this article, we present a general framework that combines the best of the discrete and continuous world: \textit{smooth} and \textit{compactly supported} basis functions, which are defined in the \textit{continuous} domain satisfying the \textit{interpolation condition} and allowing to \textit{vary} the resolution of a constructed shape. In interactive shape modeling, these properties allow for the following key attributes:

\begin{itemize}
  \item Organic shape modeling: smoothness enables a continuously-defined tangent plane and Gaussian curvature at any point on the surface, which facilitates realistic texturing and rendering of shapes; 
  \item Local shape control: compact support combined with the interpolation property of the basis functions guarantees precise and direct shape interaction and an intuitive modeling process.
  \item Detailed surfacing: few parameters are required at the initial stage of modeling, while varying the resolution of the shape allows the user to increase the number of control points when more details are to be modeled.
\end{itemize}

\noindent Our framework consists of a new family of compactly supported interpolators that are linear combinations of shifted exponential B-splines on the half-integer grid. This allows us to harness useful properties of B-splines which are then transferred to the interpolators. We first derive general results and define the construction problem together with necessary constraints and conditions. We then establish relevant reproduction properties and show that, under suitable conditions, the integer shifts of the generators form a Riesz basis, which guarantees a unique and stable representation of the parametric shapes used in practice. The generators are compactly supported. Their degree of regularity can be increased at will. {Based on extensive experimentations, we conjecture that the proposed construction always yields bona fide interpolators.}

We further propose an algorithm to change the resolution of the generators which, in turn, allows us to change the resolution of the shapes. This demands that the generators be expressed as a linear combination of finer-resolution basis functions. For this purpose, we propose a refinement scheme associated to our generators by introducing a ``pre-refinement'' step such that the resulting refinement converges to the interpolator itself. In particular, we illustrate our theory by characterizing a family of symmetric and smooth interpolators that are at least in $\mathcal C^1$ and have compact support.

Finally, we present examples of applications that involve character design (Figure~\ref{fig: fancy shapes}) as well as the design of idealized parametric shapes (Figure~\ref{fig:illustrations}).

More specifically, Sections~\ref{sec: general characterization} and~\ref{sec: interpolators for practice} are the main technical contributions, whereas in Section~\ref{sec: applications} we present practical applications which motivate this article.

\subsection{Related Work}
\noindent Recently, a method to build piecewise-polynomial interpolators has been presented in~\cite{Beccari13,Antonelli2014} and its bivariate generalization was proposed in~\cite{ANTONELLI2016}. The present work is the continuation of our previous efforts to, first, generalize the popular Catmull-Rom~\cite{Catmull74} and Keys~\cite{Keys81,Meijering2003} interpolators for practical applications~\cite{schmitter15,schmitter16b,schmitter15b,schmitter16} and, next, to go one step further and construct families of interpolators that allow varying the resolution of a shape~\cite{Badoual2015,Badoual2016}. Here, the novelty w.r.t.~\cite{schmitter16} is that the presented framework allows one to vary the resolution of shapes which facilitates shape design in practice, as illustrated in Section~\ref{sec: character design}.

\section{Review of Exponential B-Splines}\label{sec: b-splines}

\noindent We briefly review the link between exponential B-splines and differential operators. This is crucial to understand the properties of the proposed family of splines. For a more in-depth characterization of exponential B-splines, we refer the reader to~\cite{Unser0503}. 

\subsection{Notation}\label{Sec:notation}

\noindent
We describe the list of roots associated to an exponential B-spline as $\boldsymbol{\alpha}=(\alpha_1, \dots,\alpha_{n_0})$.
Likewise, we write $\alpha_n\in\boldsymbol\alpha$ to signify that one of the components of $\boldsymbol\alpha$ is $\alpha_n$. 
{The symbol $n_d$ refers to the number of distinct roots of $\boldsymbol \alpha$, which are denoted by $\alpha_{(1)},\dots,\alpha_{(n_d)}$ with the multiplicity of $\alpha_{(m)}$ being $n_{(m)}$ and ${\sum_{m=1}^{n_d}n_{(m)}=n_0}$.
The identity and derivative operators are denoted by $\mathrm{I}$ and $\mathrm{D}=\frac{\mathrm{d}}{\mathrm{d}t}$, respectively. We denote by $f(\cdot)$ a continuously defined function where the dot in parantheses represents the variable and by $a=(a[n])_{n\in \mathbb{Z}}$ a discrete sequence. The imaginary complex unit $\mathrm i$ satisfies $\mathrm{i}^2=-1$, while the Fourier integral of a function $f$ is denoted by $\widehat f(\omega)=\int_{\mathbb{R}}f(t)\mathrm{e}^{-\mathrm{i}\omega t} \mathrm{d}t$. Finally, the continuous convolutions between two functions $f$ and $g$ is  defined by $(f\ast g)(t)=\int_{\mathbb{R}}f(t-u)g(u)\mathrm{d}u$, and the discrete convolution between two sequences $a$ and $b$, is defined by $ (a \ast b)[k]=\sum\limits_{n=-\infty}^{+\infty}a[k-n]b[n]$, respectively. Furthermore, we use bold font to denote parametric shapes such as for example a 2D curve $\boldsymbol r(t)=(r_x(t),r_y(t))$.}

\subsection{{Exponential B-Spline and the Reproduction of Exponential Polynomials}} \label{subsection:operators}
\noindent The exponential B-spline with parameter $\boldsymbol \alpha$ is defined in the Fourier domain as
\begin{equation}\label{eq: B-spline fourier}
\widehat{\beta}^+_{\boldsymbol\alpha}(\omega)
=\prod\limits_{n=1}^{n_0}\frac{1-\mathrm{e}^{\alpha_n-\mathrm{i}\omega}}{\mathrm{i}\omega-\alpha_n}.
\end{equation}

\noindent The function $\beta^+_{\boldsymbol{\alpha}}$ is compactly supported with support $[0,n_0]$~\cite[Section III-A]{Unser0503}. 
%
%

\noindent We denote by $\beta_{\bm{\alpha}}$ the corresponding \textit{centered} (hence, non-causal) exponential B-spline, whose support is $[-n_0/2, n_0/2]$. We have therefore
\begin{equation}\label{eq:betacausalcentered}
 \beta_{\bm{\alpha}} (t) = \beta^+_{\bm{\alpha}} (t + n_0 /2),
 \end{equation}
with $\beta^+_{\bm{\alpha}}$ the \textit{causal} B-spline defined in~\eqref{eq: B-spline fourier}. 
The reason for introducing centered B-splines is that we shall define interpolators that are symmetric around the origin and, hence, centered. 

{It is well known that the exponential B-spline $\beta_{\bm{\alpha}}$ is intimately linked to the differential operator 
\begin{equation}\label{eq: operator la}
	\mathrm{L}_{\bm{\alpha}} = (\mathrm{D} - \alpha_1 \mathrm{I}) \ldots (\mathrm{D} - \alpha_{n_0} \mathrm{I} ),
\end{equation}
which implies that $\beta_{\boldsymbol\alpha}$ is able to reproduce the functions $p_0$ in the null space of $L_{\alpha}$ defined as $\mathrm{L}_{\alpha} p_0 = 0$.
As a consequence, exponential B-splines can reproduce exponential polynomials that live in the space~\cite[Section III-C-2]{Unser0503}}
\begin{equation}\label{eq: reproducing space}
{\mathrm{span}\{t^{n-1}\mathrm{e}^{\alpha_{(m)}t}\}_{m=1,\dots,n_d; \ n=1,\dots,n_{(m)}}.}
\end{equation}

\section{General Characterization of the Interpolator}\label{sec: general characterization}
\noindent We consider generators that are constructed as a sum of half-integer shifted versions of a given exponential B-spline $\beta_{\bm{\alpha}}$. 
\begin{definition}\label{def: basis function}
For a sequence $\lambda \in \ell_1(\Z)$ and $\bm{\alpha}$ a vector of roots, we define
\begin{equation} \label{eq:phialpha}
	\phi_{\lambda,\bm\alpha}(t):= \sum_{n\in \Z} \lambda[n] \beta_{\bm{\alpha}} \left( t -  \frac{n}{2}\right).
\end{equation}
\end{definition}

\noindent In the frequency domain, we then have
\begin{equation}\label{eq: interpolator fourier}
	\widehat{\phi}_{\lambda, \bm{\alpha}} (\omega) =  \left( \sum_{n\in \mathbb{Z}} \lambda[n] \mathrm{e}^{ - \mathrm{i} \omega n /2} \right) \widehat{\beta}_{\bm{\alpha}} (\omega). 
\end{equation}

\noindent In what follows, we state the desired mathematical properties that the generator $\phi_{\lambda,\bm\alpha}$ should satisfy. 
\begin{enumerate}[I]
	\item 
	The generator $\phi_{\lambda,\bm\alpha}$ is interpolatory in the sense that, for any function $f \in \mathrm{span}\{ \phi_{\lambda,\bm\alpha} ( \cdot - k )\}_{k \in \Z }$,
we have $f(t) = \sum_{k\in \Z} f[k] \phi_{\lambda,\bm\alpha}(t-k)$.
	This is equivalent to the interpolation condition 
	\begin{equation}
	 \phi_{\lambda,\bm\alpha}(t)\big|_{t=k}=\delta[k]=\left\{
      \begin{aligned}
        1 & \mbox{ if } k=0\\
        0 & \mbox{, otherwise}\\
      \end{aligned}
    \right.,
	\end{equation}
	where $\delta[k]$ represents the Kronecker delta.	
	\item
	The generator $\phi_{\lambda,\bm\alpha}$ is compactly supported, which implies that the sequence $\lambda$ has a finite number of non-zero values.
	\item
	The function $\phi_{\lambda,\bm\alpha}$ is smooth with at least a continuous derivative.
	\item
	The family of the integer shifts of the generator $\{\phi_{\lambda,\bm\alpha}(\cdot-k)\}_{k\in\Z}$ forms a Riesz basis.
	\item
	The generator $\phi_{\lambda,\bm\alpha}$ preserves the reproduction properties of the associated exponential B-spline $\beta_{\bm{\alpha}}$ in the sense that it is capable of reproducing the exponential polynomials in the null-space of the operator $\mathrm L_{\boldsymbol{\alpha}}$ defined in~\eqref{eq: operator la}.
	\item
	The generator $\phi_{\lambda,\bm\alpha}$ allows one to represent shapes at various resolutions.
\end{enumerate}

\noindent We choose equispaced half-integer shifts of the exponential B-splines in Definition~\ref{def: basis function}. The reason is that our problem has no solution using only integer shifts under Conditions I), II), and III): There is no smooth and compactly supported interpolator of the form $\sum_{k \in \Z} \lambda[k] \beta_{\bm\alpha} ( t - k)$. This can easily be verified; for example, by plugging any polynomial B-spline into Definition~\ref{def: basis function} and using integer shifts while imposing the interpolation conditions: It turns out that there are not enough degrees of freedom to solve the problem due to the compact support of the B-splines as well as the smoothness condition, which forces the degree of the B-spline to be greater than $1$. Furthermore, by using half-integer shifts, we guarantee that our solution lives in the spline space of the next finer resolution; a property that can be exploited in practice, as detailed in Section~\ref{sec: varying resolution}.

\subsection{Riesz Basis}
\noindent We consider the space 
\begin{equation}
	V(\phi_{\lambda,\bm\alpha} )  = \left\{ \sum_{n\in \Z} c[n] \phi_{\lambda,\bm{\alpha}}(\cdot - n), \   c \in \ell_2(\Z)  \right\}
\end{equation}
of functions that is generated by the integer shifts of $\phi_{\lambda,\bm{\alpha}}$. Our requirement is that the family of functions $\{ \phi_{\lambda,\bm{\alpha}}(\cdot - n) \}_{n\in \Z}$ forms a Riesz basis of $V(\phi_{\lambda,\bm{\alpha}})$, which ensures that the representation of a function in $V(\phi_{\lambda,\bm{\alpha}})$ is stable and unique. We show in this section that this is the case if $\{ \beta_{\bm{\alpha}}(\cdot - n) \}_{n\in \Z}$ is itself a Riesz basis and if $\phi_{\lambda,\bm{\alpha}}$ is interpolatory.

\begin{definition}
	The family $\{ \varphi_n \}_{n\in \mathbb{Z}}$ of functions forms a Riesz basis if
	\begin{equation} \label{eq:rieszbasiscondition}
		A \lVert c \rVert_{\ell_2(\Z)}  \leq  \left\lVert  \sum_{n\in \Z} c[n] \varphi_n \right\rVert_{L_2(\R)} \leq B \lVert c \rVert_{\ell_2(\Z)}
	\end{equation}
	for some constants $A,B>0$ and any sequence $c = (c[n])_{n\in \Z} \in \ell_2(\Z)$.
\end{definition}

When $\varphi_n = \varphi(\cdot - n)$,~\eqref{eq:rieszbasiscondition} is equivalent to the Fourier-domain condition
	\begin{equation} \label{eq:rieszphi}
		A^2 \leq \sum_{k\in \mathbb{Z}} \lvert \widehat{\varphi}(\omega - 2 k \pi) \rvert^2 \leq B^2
	\end{equation}
	for any $\omega\in \R$~\cite{Unser2000}. The family $\{ \beta_{\bm\alpha} ( \cdot - n )\}_{n\in \Z}$ is a Riesz basis when $\bm\alpha$ is such that $\alpha_n  - \alpha_m \neq 2 k \pi \mathrm{i}$, $k\in \mathbb{Z}$, for any pair of distinct purely imaginary roots $\alpha_m, \alpha_n \in \bm\alpha$~\cite[Theorem 1]{Unser0503}. 


\begin{proposition} \label{prop:riesz}
	Let $\bm\alpha$ be such that $\alpha_n  - \alpha_m \neq 2 k \pi \mathrm{i}$, $k\in \mathbb{Z}$, for any pair of distinct purely imaginary roots $\alpha_m,  \alpha_n \in \bm\alpha$. 	
	For any sequence $\lambda \in \ell_1(\Z)$, if the basis function $\phi_{\lambda, \bm{\alpha}}$ is interpolatory, then the family $\{ \phi_{\lambda,\bm{\alpha}} (\cdot - n) \}_{n\in \mathbb{Z}}$ is a Riesz basis.
\end{proposition}
 
\noindent The proof is given in~\ref{prop:appendixriesz} as well as an estimate of the Riesz Bounds. 

\subsection{Reproduction Properties}
\begin{proposition}\label{prop: lambda reproduction}
Let $\bm\alpha$ be a vector of roots. 
We assume that $ \lambda \in \ell_1(\Z)$ satisfies the conditions
\begin{align} \label{eq:conditionslambdareproductionexistence}
	\sum_{n\in \Z} \lvert \lambda[n]\rvert  \mathrm{e}^{- \alpha n /2} &<  \infty,  \\ 
	\sum_{n\in \Z} \lambda[n] \mathrm{e}^{- \alpha n /2} &\neq  0 \label{eq:conditionslambdareproduction}
\end{align}
for every $\alpha \in \bm\alpha$. Then, the basis function $\phi_{\lambda,\bm{\alpha}}$ has the same reproduction properties as the corresponding exponential B-spline $\beta_{\bm\alpha}$. In particular, it reproduces the exponential polynomials 
\begin{equation}
t^{n-1}\mathrm{e}^{\alpha_{(m)}t}
\end{equation}
for $m=1,\dots,n_d$ and $n=1,\dots,n_{(m)}$, with the notations of Section~\ref{Sec:notation}.
\end{proposition}

Note that~\eqref{eq:conditionslambdareproductionexistence} is always satisfied as soon as $\phi_{\lambda,\boldsymbol{\alpha}}$ is compactly supported. The proof of Proposition~\ref{prop: lambda reproduction} is given in~\ref{Proof of Prop two}.

\subsection{Regularity} \label{sec:regularity}
\noindent From Definition~\ref{def: basis function}, it immediately follows that $\phi_{\lambda, \bm\alpha}$ has the same regularity as the exponential B-spline $\beta_{\bm\alpha}$ if $\lambda \neq 0$.
Hence, $\phi_{\lambda,\bm\alpha}$ belongs to $\mathcal{C}^{n_0-2}$~\cite[Section III-A]{Unser0503}.

\subsection{Varying the Resolution of the Generator}\label{sec: varying resolution}
\noindent The causal exponential B-spline $\beta_{\bm \alpha}^+$ is refinable, in the sense that its dilation by an integer $m$ can be expressed as a linear combination of $\beta_{\bm\alpha/m}^+(\cdot - k)$. This is what we refer to as the \textit{resolution} of the basis function. We shall see how this property translates to the function $\phi_{\lambda,\bm\alpha}$.
For this purpose, we first revisit the $m$-scale relation for exponential B-splines. For convenience, we express the corresponding terms with respect to causal (non-centered) B-splines. 
In practice, we always consider symmetric interpolators $\phi_{\lambda, \bm{\alpha}}$ with support $[-(n_0 -1), n_0 - 1]$  (see Section~\ref{sec: interpolators for practice}). Therefore, we define the shifted and causal version of the interpolator as
\begin{equation} \label{eq:phi_causal}
\phi_{\lambda,\bm\alpha}^+ (t) = \phi_{\lambda,\bm\alpha}( t - (n_0-1)).
\end{equation}
Every causal formula is easily adapted to the centered case by applying a shift similar to~\eqref{eq:phi_causal}.
We follow the notations of~\cite{Unser0503}, where an in-depth discussion on the refinability of exponential B-splines can be found.
 
As shown in~\cite[Section IV-D]{Unser0503}, the dilation by an integer $m\in \mathbb{N}\setminus\{0\}$ of an exponential B-spline is expressed in the space domain as
\begin{equation}\label{eq: dilation relation}
\beta^+_{\boldsymbol\alpha} \left(\frac{t}{m}\right)=\sum_{k\in\mathbb{Z}}h_{\frac{\boldsymbol\alpha}{m},m}[k]\beta^+_{\frac{\boldsymbol\alpha}{m}}(t-k),
\end{equation}
where the \emph{refinement filter} $h_{\bm \alpha, m}$ is specified by its Fourier transform as
\begin{equation}
\label{eq: refinement mask}
H_{\boldsymbol\alpha,m}(\mathrm{e}^{\mathrm{i}\omega})=\frac{1}{m^{n_0-1}}\prod_{n=1}^{n_0}\bigg(\sum_{k=0}^{m-1}\mathrm{e}^{\alpha_nk}\mathrm{e}^{- \mathrm{i} k \omega}\bigg).
\end{equation}

\noindent As we shall see, it is impossible to establish a similar relation for the interpolator $\phi^+_{\lambda,\bm\alpha}$. However, we can exploit the refinability of the corresponding spline $\beta_{\bm\alpha}^+$ to express the dilation of $\phi^+_{\lambda,\bm\alpha}$. 

For $\bm\alpha$ a vector of roots, $\lambda \in \ell_1(\Z)$, and $m_0$ an \textit{even} integer, we define the digital \emph{pre-filter} $g_{\lambda,\bm\alpha,m_0}$ by its Fourier transform
\begin{equation}
	G_{\lambda,\bm\alpha,m_0} (\mathrm{e}^{\mathrm{i}\omega}) =   \mathrm{e}^{- \mathrm{i}\omega m_0 (n_0 / 2 - 1)}  \left( \sum_{n \in \Z} \lambda[n] \mathrm{e}^{- \mathrm{i}\omega n m_0 / 2} \right) H_{\frac{\bm\alpha}{m_0},m_0}(\mathrm{e}^{\mathrm{i}\omega}).
\end{equation}
The term $ \mathrm{e}^{- \mathrm{i}\omega m_0 (n_0 / 2 - 1)}$ is due to the fact that $\beta_{\bm\alpha}$ and $\phi_{\lambda,\bm\alpha}$ do not have the same support in general. The pre-filter allows us to express $\phi^+_{\lambda,\bm\alpha}$ dilated by $m_0$ as a linear combination of the refined shifted B-splines $\beta_{\frac{\bm\alpha}{m_0}}^+(\cdot - k)$. 
Note that $G_{\lambda,\bm\alpha,m_0}$ is a valid Fourier transform of a digital filter (\emph{i.e.}, a function of $\mathrm{e}^{\mathrm{i} \omega}$) only for even $m_0$.

\begin{proposition} \label{prop:pre filtering}
	Let $\bm\alpha$ be a vector of roots, $\lambda \in \ell_1(\Z)$, and $m_0$ be an even integer. 
	Then, we have
	\begin{equation}\label{eq: dilation relation for phi}
	\phi^+_{\lambda, \boldsymbol\alpha} \left(\frac{t}{m_0}\right)=\sum_{k\in\mathbb{Z}}g_{\lambda, \frac{\boldsymbol\alpha}{m_0},m_0}[k]\beta^+_{\frac{\boldsymbol\alpha}{m_0}}(t-k).
	\end{equation}
	\end{proposition}

The proof is given in~\ref{Proof of Prop three}.

\subsubsection{Modified Refinement Scheme Based on Exponential B-Splines}
\noindent Using Proposition~\ref{prop:pre filtering}, we are able to express a function which is constructed with the interpolator $\phi^+_{\lambda, \boldsymbol\alpha}$ in an exponential B-spline basis. Starting with the samples $c[k]=f(t)|_{t=k\in\mathbb{Z}}$ of a continuously defined function $f(\cdot)$ that can be \textit{perfectly} reconstructed, \textit{i.e.}, $f\in\mathrm{span}\{\phi^+_{\lambda,\boldsymbol\alpha}(\cdot-k)\}_{k\in\mathbb Z}$, we have
	\begin{align}\label{eq: change of basis}
		f(t)&=\sum_{k\in\mathbb{Z}}{c}[k]\phi^+_{\lambda,\boldsymbol\alpha}(t-k)  \nonumber \\
&=\sum_{k\in\mathbb{Z}}{c}[k]
\sum_{l\in\mathbb{Z}}
g_{\lambda,\frac{\boldsymbol\alpha}{m_0},m_0}[l]\beta^+_{\frac{\boldsymbol{\alpha}}{m_0}}(m_0(t-k)-l) \nonumber  \\
&=\sum_{k\in\mathbb{Z}}\sum_{l\in\mathbb{Z}}{c}[k]g_{\lambda, \frac{\boldsymbol\alpha}{m_0},m_0}[l]\beta^+_{\frac{\boldsymbol{\alpha}}{m_0}}(m_0t-m_0k-l)  \nonumber  \\
&=\sum_{l\in\mathbb{Z}}c_0 [l]\beta^+_{\frac{\boldsymbol{\alpha}}{m_0}}(m_0t-l)
	\end{align}		
with 
	\begin{equation}\label{eq: modified subdiv}
		c_0[l]=\bigg({c}_{\uparrow m_0}\ast g_{\lambda,\frac{\boldsymbol\alpha}{m_0},m_0}\bigg)[l],
	\end{equation}
		where $\uparrow m_0$ denotes upsampling by a factor $m_0$ defined as
\begin{equation}\label{eq: upsampling}
{c_{\uparrow_{m_0}}[k]=\left\{
\begin{array}{ll}
c[n] \text{, }k=m_0 n\\
0\text{, otherwise.}
\end{array}
\right.}
\end{equation}

Equation~\eqref{eq: change of basis} shows that a function that is originally expressed in the basis generated by $\phi^+_{\lambda, \boldsymbol\alpha}$ can be expressed in a corresponding exponential B-spline basis with respect to a finer grid. This suggests that, after having performed the change of basis described by~\eqref{eq: change of basis}, the resolution of $f$ can be further refined by applying the standard iterative B-spline refinement rules. At this point, it is interesting to take a deeper look into the relation between the interpolated function $f$ and the sequence $c$ of samples as we iteratively refine it. As will become apparent in the application-oriented Section~\ref{sec: applications}, a parametric shape is described by coordinate functions whose samples build 2D or 3D vectors of \textit{control points}. Repositioning of these control points allows us to locally modify the shape, while the iterative refinement of the control points allows us to iteratively increase the local control over the shape. Hence, for practical purposes, it is convenient to study the convergence of the refinement process as the number of iterations becomes large. Proposition~\ref{prop: refinement algo} describes the refinement scheme and provides the corresponding convergence result.

\begin{proposition}\label{prop: refinement algo}
Let $\bm\alpha$ be a vector of roots and $\lambda \in \ell_1(\Z)$. For a continuous function $f$ with samples $f(t)|_{t=k\in\mathbb{Z}}=c[k]$ and the integers $m,m_0$, with $m_0$ being even, we consider the iterative scheme specified by
	\begin{enumerate}
		\item \textit{pre-filter step:} $c_0[k] = (g_{\lambda, \frac{\bm\alpha}{m_0}, m_0} \ast {c}_{\uparrow m_0}) [k]$;
		\item \textit{iterative steps:} for $n\geq 1$, $c_n[k] = (h_{\frac{\bm\alpha}{m_0m^n}, m} \ast {(c_{n-1})}_{\uparrow m}) [k]$,
	\end{enumerate}
{where $\uparrow m$ denotes upsampling by a factor $m$ as defined in~\eqref{eq: upsampling}.}
	Then, the iterative scheme is convergent, in the sense that 
	\begin{equation}\label{eq: convergence result}
		\lim_{n\rightarrow\infty}\sum_{k \in \Z} c_n[k] \delta\left( m_0m^n t - k\right) = f(t),
	\end{equation}
	{where $\delta$ is the Dirac distribution.}
\end{proposition}

The proof is given in Appendix~\ref{Proof of Prop four}.

\subsubsection{Example}
\begin{figure}[htb]
\centering
\includegraphics[width=8cm]{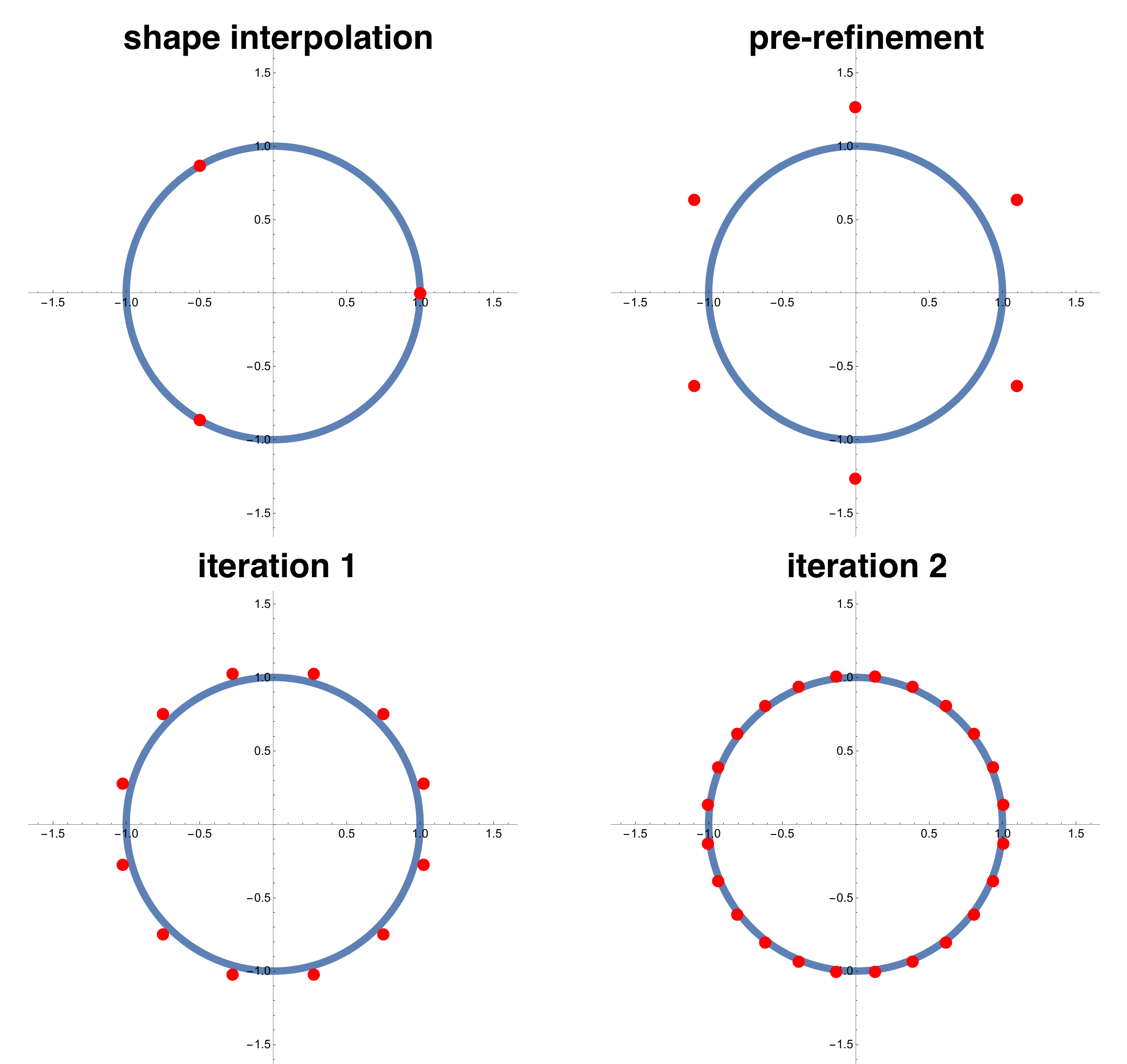}
\caption{Refined circle. The parametric circle is first constructed using the proposed interpolator and $\boldsymbol\alpha=(0,\frac{2\mathrm{i}\pi}{3},-\frac{2\mathrm{i}\pi}{3})$ (top left). At the first iteration, the ``pre-refinement'' mask is applied to the initial control points according to~\eqref{eq: modified subdiv} (top right), whereas at the subsequent iterations the standard refinement mask for exponential B-splines~\eqref{eq: refinement mask} is applied (bottom, from left to right). In the bottom right, we see how the iterative process converges towards the continuously defined circle.}
\label{fig: subdiv circle}
\end{figure}

\noindent We illustrate how to refine the resolution of a circular pattern by applying Proposition~\ref{prop: refinement algo}. To efficiently take advantage of the interpolation property, we apply the ``pre-refinement" step~\eqref{eq: modified subdiv} at the first iteration. For the subsequent iterations, we apply the standard refinement given by~\eqref{eq: refinement mask} as described by Proposition~\ref{prop: refinement algo}. By doing so, we see that the iterative scheme converges towards the circle $\boldsymbol r(t) = \sum_{l\in\mathbb{Z}}{\boldsymbol c}_{0}[l]\beta^+_{\frac{\boldsymbol{\alpha}}{m_0}}(m_0 t-l) = \sum_{k\in\mathbb{Z}}{\boldsymbol r}[k]\phi^+_{\boldsymbol\alpha}(t-k)$. The result of the algorithm is shown in Figure~\ref{fig: subdiv circle}. In~\ref{sec: repro ellipse}, we provide the details on how to reconstruct the circle with our framework.

\bigskip
 
 \newpage
 \clearpage
 
\section{Construction of a Family of Compactly Supported Interpolators in Practice}\label{sec: interpolators for practice}
\noindent It is known that there exists no exponential B-spline $\beta_{\bm{\alpha}}$ that is interpolatory and smooth (\textit{i.e.}, at least in $\mathcal C^1$) at the same time. Our goal here is to construct a compactly supported generator function that has the same smoothness and reproduction properties as $\beta_{\bm{\alpha}}$, while also being interpolatory. In order to meet the smoothness constraints, we require the number of elements of $\bm\alpha$ to be $n_0 \geq 3$ in accordance with the construction detailed in Section~\ref{sec:regularity}. Furthermore, we want the interpolator to be real-valued and symmetric, which implies that the elements of $\bm\alpha$ are either zero or come in complex conjugate pairs~\cite{Unser0503}. Using Definition~\ref{def: basis function} and the conditions described in Section~\ref{sec: general characterization}, we are looking for the interpolator with {minimal} support.

\subsection{Introductory Example: The Quadratic B-Spline} \label{subsec:quadraticBspline} 
\noindent We illustrate the concept with a simple example that uses quadratic polynomial B-splines, which are constructed with $\boldsymbol\alpha=\boldsymbol\alpha_0=(0,0,0)$ in~\eqref{eq: B-spline fourier} and whose support is of size 3. The interpolation constraint combined with the half-integer shifts demand that $\lambda$ contains at least three non-zero values to have enough degrees of freedom. This also implies that the compactly-supported interpolator is constructed with no more than three non-zero elements of $\lambda$. {Moreover, since the solution that fulfills the conditions stated in Section~\ref{sec: general characterization} is unique, the interpolator is of minimal support}. To satisfy the symmetry constraints, we center the shifted B-splines around the origin and enforce $\lambda[1]=\lambda[-1]$. Hence, our generator must take the form
\begin{equation}
\begin{split}
\phi_{\lambda,\bm{\alpha_0}}(t)&=\lambda[1]\beta_{\bm{\alpha_0}}(t-\frac{1}{2})+\lambda[0]\beta_{\bm{\alpha_0}}(t)+\lambda[-1]\beta_{\bm{\alpha_0}}(t+\frac{1}{2})\\
&=\lambda[0]\beta_{\bm{\alpha_0}}(t)+\lambda[1]\big(\beta_{\bm{\alpha_0}}(t-\frac{1}{2})+\beta_{\bm{\alpha_0}}(t+\frac{1}{2})\big).
\end{split}
\end{equation}

\noindent Since $\boldsymbol\alpha_0$ has $n_0=3$ elements, the support of the interpolator is $N = 2(n_0 - 1)=4$. The interpolator itself is supported in $[-(n_0 - 1),(n_0 - 1)]=[-2,2]$. 
The interpolation condition is expressed as 
$
\begin{cases}
\phi_{\lambda,\bm{\alpha_0}}(0)=1\\
\phi_{\lambda,\bm{\alpha_0}}(1)=0
\end{cases}
$.
We define the matrix
\begin{equation*}
\mathbf{A}_{\boldsymbol\alpha_0}=
\begin{pmatrix}
 \beta_{\bm{\alpha_0}}(0) &  \beta_{\bm{\alpha_0}}(-1/2)+ \beta_{\bm{\alpha_0}}(1/2)\\
  \beta_{\bm{\alpha_0}}(1) &  \beta_{\bm{\alpha_0}}(1-1/2) +  \beta_{\bm{\alpha_0}}(1+1/2)
\end{pmatrix}
=
\left(
\begin{array}{cc}
 \frac{3}{4} & 1 \\
 \frac{1}{8} & \frac{1}{2} \\
\end{array}
\right)
\end{equation*}
and rewrite the interpolation constraint as 
$(\lambda[0],\lambda[1]) = \bm{\mathrm{A}}_{\bm{\alpha_0}}^{-1} (1,0)=(1,-\frac{1}{2})$. The resulting interpolator is shown in Figure~\ref{fig: interpolators}.

\subsection{The General Case}
\noindent In what follows, we only consider vectors of poles $\bm\alpha$ for which $\alpha_n  - \alpha_m \neq 2 k \pi \mathrm{i}$, $k\in \mathbb{Z}$ for all pairs of distinct, purely imaginary roots $\alpha_m, \alpha_n \in \bm\alpha$ (Riesz Basis property). 
We generalize the above example to construct symmetric and compactly supported interpolators of any order and that are of the form
\begin{equation} \label{eq:varphialpha}
	\phi_{\lambda,\bm\alpha}(t):= \lambda[0] \beta_{\bm{\alpha}} ( t )+\sum_{n=1}^{n_0-2} \lambda[n] \left( \beta_{\bm{\alpha}} ( t - n/2) +  \beta_{\bm{\alpha}} ( t + n/2)\right),
\end{equation}
whose support is included\footnote{The support is exactly $[- N/2, N/2] $ when $\lambda[n]$ is non-zero for $n=0,\ldots, (n_0-2)$, which is always the case in the examples we have considered.} in $[- N/2, N/2] = [-(n_0 - 1), n_0 -1]$. We easily pass from the general representation~\eqref{eq:phialpha} to~\eqref{eq:varphialpha}, adapted to the symmetric and compactly supported case, by setting $\lambda[n] = 0$ when $\lvert n \rvert \geq n_0 -1$ (support condition) and $\lambda[-n] = \lambda[n]$ for every $n$ (symmetry condition).

The function $\phi_{\lambda,\bm\alpha}$ is interpolatory if and only if 
\begin{equation} \label{eq:systemlambda}
	\phi_{\lambda,\bm\alpha}(0) = 1  \text{ and }
	\phi_{\lambda,\bm\alpha}(1) = \cdots = \phi_{\lambda,\bm\alpha}(n_0-2) =0.
\end{equation}
This defines a linear system with $(n_0 -1)$ unknown non-zero elements of $\lambda$, $\{\lambda[0], \ldots , \lambda[n_0 -2]\}$, and $(n_0 -1)$ equations. 	The system~\eqref{eq:systemlambda} has a solution if the matrix $\bm{\mathrm{A}}_{\bm{\alpha}} \in \R^{(n_0-1)\times(n_0-1)}$ defined for $k,l = 0, \ldots , (n_0-2)$ by
	\begin{equation}	 \label{eq:Aalpha}
		[\bm{\mathrm{A}}_{\bm{\alpha}}]_{k+1,l+1} =
		\begin{cases}
		 \beta_{\bm{\alpha}} ( k ) & \quad\text{if}\quad l=0  \\
		 \beta_{\bm{\alpha}} ( k - l/2) + \beta_{\bm{\alpha}} ( k + l/2) &	\quad\text{else}\quad	
		 \end{cases}
	\end{equation}
	is invertible. In this case, we have 
	\begin{equation} \label{eq:conditionlambda}
		\lambda = (\lambda[0],\ldots,\lambda[n_0-2]) = \bm{\mathrm{A}}_{\bm{\alpha}}^{-1} (1,0,\ldots,0).
	\end{equation}
	
\noindent Knowing $\bm{\alpha}$, we can easily check if the matrix $\bm{\mathrm{A}}_{\bm{\alpha}}$ is invertible, which is the case for all the examples that we tested (we have already seen that it is true for $\bm\alpha = (0,0,0)$ is Section~\ref{subsec:quadraticBspline}). From~\eqref{eq:conditionlambda}, we see that $\lambda$ is completely determined by $\bm\alpha$. This motivates Definition~\ref{def: interpolator time domain}.
	
\begin{definition} \label{def: interpolator time domain}
Let $\bm\alpha$ be a vector of roots whose elements are either zero or come in pairs with opposite signs. If the matrix $\bm{\mathrm{A}}_{\bm{\alpha}}$ defined in~\eqref{eq:Aalpha} is invertible, then the interpolatory basis function $\phi_{\bm\alpha}$ is defined as
		\begin{equation}
			\phi_{\bm\alpha} := \phi_{\lambda,\bm\alpha},
		\end{equation}
		with $\lambda$ defined by~\eqref{eq:conditionlambda}.  
	\end{definition}

We conjecture that the matrix $\bm{\mathrm{A}}_{\bm{\alpha}}$ is always invertible, and that we always can define an interpolator $\phi_{\bm{\alpha}}$ for any list of roots $\boldsymbol\alpha$. 
In the remaining of this article, we assume that $\bm{\mathrm{A}}_{\bm{\alpha}}$ is invertible and, therefore, that $\phi_{\bm\alpha}$ is well-defined. 
Under this assumption, the unicity of the vector $\lambda$ ensures that the interpolator $\phi_{\bm\alpha}$ in Definition~\ref{def: interpolator time domain} has minimal support among the interpolators of the form~\eqref{eq:phialpha}. 

\begin{figure}[htb]
\begin{centering}
\includegraphics[scale=0.4]{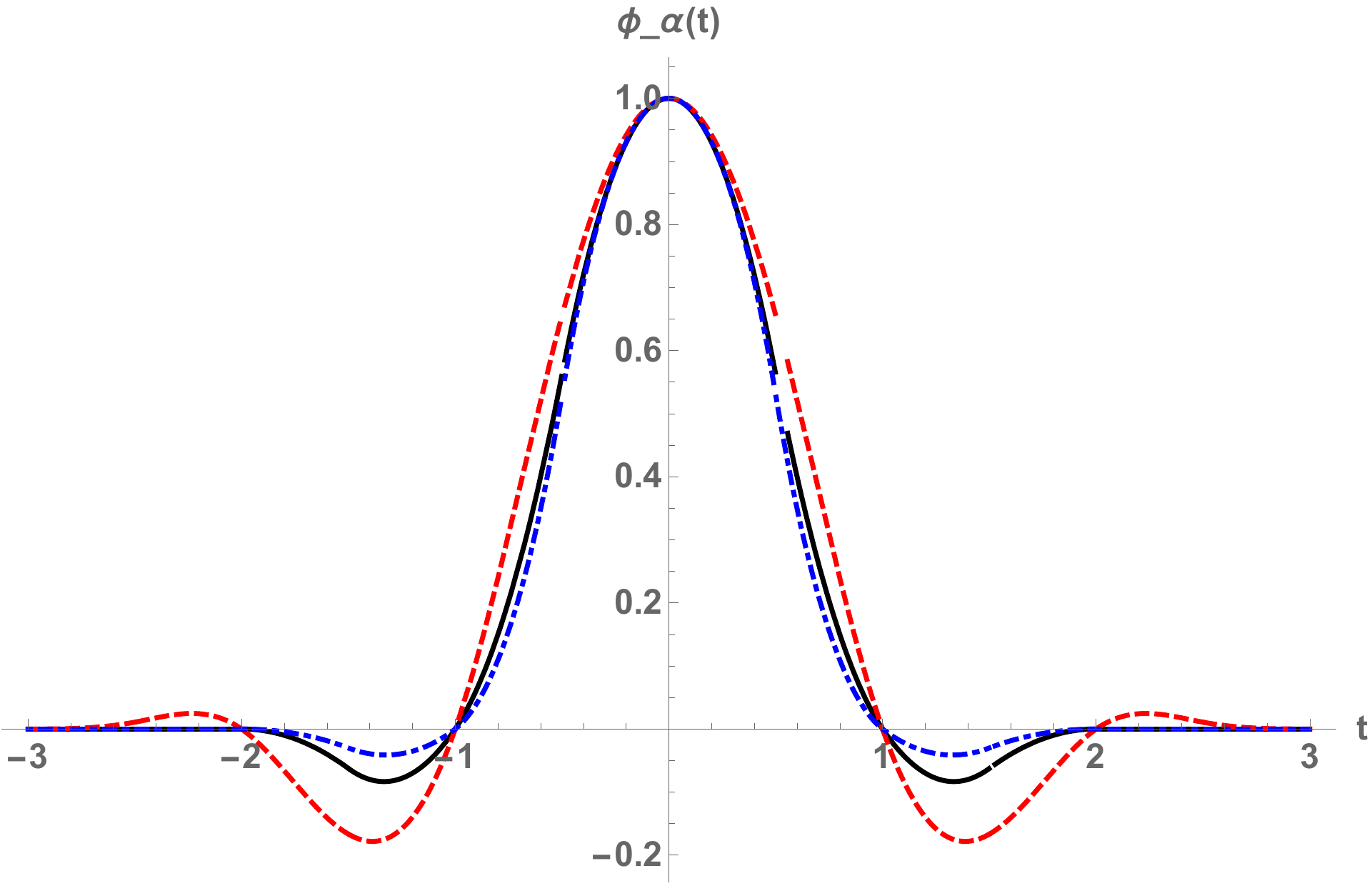}
\caption{Different types of interpolators: polynomial interpolator (black, solid curve) with $\boldsymbol\alpha=(0,0,0)$. The number of poles is equal to $3$. Trigonometric interpolator (red, dashed curve): the non-zero poles are purely imaginary and come in pairs with opposite signs (\textit{e.g.}, $\boldsymbol\alpha = (0,0, \frac{\mathrm{i}2\pi}{3},-\frac{\mathrm{i}2\pi}{3})$). Hyperbolic interpolator (blue, dot-dashed curve): the non-zero poles are real and come in pairs with opposite signs (\textit{e.g.}, $\boldsymbol\alpha = (0, \frac{2\pi}{3},-\frac{2\pi}{3})$).}
\label{fig: interpolators}
\end{centering}
\end{figure}

In practice, the type of interpolator that needs to be constructed depends on the parametric shape that is represented. For instance, for a rectangular surface, a polynomial interpolator is required and the vector $\boldsymbol\alpha$ of roots will have to consist of zeros. If instead we aim at representing circles, spheres, or ellipsoids (see Section~\ref{sec: applications}), whose coordinate functions are trigonometric, we need to construct
interpolators that preserve sinusoids. Therefore, $\boldsymbol\alpha$ will contain pairs of purely imaginary roots with opposite signs. Similarly, we can reproduce hyperbolic shapes by picking an $\boldsymbol\alpha$ that contains pairs of real roots with opposite signs. If an interpolator is required to reproduce both trigonometric and polynomial shapes, \textit{e.g.}, to construct a cylinder, then the corresponding polynomial and trigonometric root vectors are concatenated to construct $\boldsymbol\alpha$. Examples of different interpolators are shown in Figure~\ref{fig: interpolators}.

\noindent We now summarize the properties of the generator $\phi_{\bm\alpha}$ for $\bm\alpha$ a vector of roots of size $n_0 \geq 3$ such that $\alpha_n  - \alpha_m \neq 2 k \pi \mathrm{i}$, $k\in \mathbb{Z}$, for any pair of distinct purely imaginary roots $\alpha_m, \alpha_n \in \bm\alpha$ .
These properties are in accordance with Conditions I to VI in Section~\ref{sec: general characterization}. 
\begin{itemize}
	\item The function $\phi_{\bm\alpha}$ is interpolatory.
	\item The function $\phi_{\bm\alpha}$ is compactly supported in $[-(n_0-1), n_0-1]$.
	\item The function $\phi_{\bm\alpha}$ has the minimal support among the interpolators that are linear combinations of shifted exponential B-splines on the half-integer grid. 
	\item The function $\phi_{\bm\alpha}$ is in $\mathcal{C}^{n_0-2}$ and therefore, at least in $\mathcal{C}^1$.
	\item The family $\{ \phi_{\bm\alpha}( \cdot - n)\}_{n\in\Z}$ is a Riesz basis.
	\item The family $\{ \phi_{\bm\alpha}( \cdot - n)\}_{n\in\Z}$ reproduces the exponential polynomials given by~\eqref{eq: reproducing space}.
	\item The function $\phi_{\bm\alpha}$ is refinable in the sense explained in Section~\ref{sec: varying resolution}. 
\end{itemize}

\noindent \textbf{Remark.} The presented interpolators are not (entirely) positive (see Figure~\ref{fig: interpolators}) and thus, do not satisfy the convex-hull-property. However, the popularity of the Catmull-Rom splines~\cite{Catmull74} in computer graphics shows that in interactive shape modeling, one prefers to use interpolators at the expense of the convex-hull property.

\section{Applications}\label{sec: applications}
\noindent In this section, we show how parametric curves and surfaces are constructed using the proposed spline bases. Such shapes can be constructed independently of the number of control points. This makes them particularly useful for deformable models where, starting from an initial configuration, one aims at approximating a target shape with arbitrary precision~\cite{Unser2000}.

\subsection{Reproduction of Idealized Shapes}
\noindent We consider curves and surfaces that are described by the coordinate functions $r_x(t)$, $r_y(t)$, and $r_z(t)$, with $t\in\mathbb R$. The coordinate functions are expressed by a linear combination of weighted integer shifts of the generator $\phi_{\boldsymbol\alpha}$. Due to the interpolation property of the generator, the weights simply correspond to the samples of the coordinate functions. Such a parametric curve is expressed as
\begin{equation}\label{eq: parametric curve}
\boldsymbol r(t)
=\begin{pmatrix}
r_x(t)\\
r_y(t)\\
r_z(t)
\end{pmatrix}
=\sum_{k\in\mathbb{Z}}\boldsymbol r[k]\phi_{\boldsymbol\alpha}(t-k),
\end{equation}
\noindent where the coefficients $\boldsymbol r[k]=(r_x[k],r_y[k],r_z[k])$ with $k\in\mathbb Z$ are the \textit{control points}. The curve~\eqref{eq: parametric curve} can be locally modified by changing the position of a single control point. The shapes that $\boldsymbol r$ can adopt (\textit{e.g.}, polynomial, circular, elliptic) depend on the properties of the generator.

One can also extend the curve model~\eqref{eq: parametric curve} to represent separable tensor-product surfaces. In this case, a surface $\boldsymbol\sigma$ is parameterized by $u,v\in\mathbb R$ as
\begin{equation}\label{eq: parametric surface}
\begin{split}
\boldsymbol\sigma(u,v)
&=
\begin{pmatrix}
\boldsymbol\sigma_x(u,v)\\
\boldsymbol\sigma_y(u,v)\\
\boldsymbol\sigma_z(u,v)
\end{pmatrix}
=\begin{pmatrix}
r_{1,x}(u)\cdot r_{2,x}(v)\\
r_{1,y}(u)\cdot r_{2,y}(v)\\
r_{1,z}(u)\cdot r_{2,z}(v)\\
\end{pmatrix}\\
&=\sum_{k\in\mathbb{Z}}\boldsymbol r_1[k]\phi_{\boldsymbol\alpha_1}(u-k)\times\sum_{l\in\mathbb{Z}}\boldsymbol r_2[l]\phi_{\boldsymbol\alpha_2}(v-l)\\
&=\sum_{k\in\mathbb{Z}}\sum_{l\in\mathbb{Z}}\underbrace{\boldsymbol r_1[k]\times\boldsymbol r_2[l]}_{\boldsymbol\sigma[k,l]}\phi_{\boldsymbol\alpha_1}(u-k)\phi_{\boldsymbol\alpha_2}(v-l),
\end{split}
\end{equation}
\noindent where ``$\times$'' denotes the element-wise multiplication of two vectors. Finally, one generalizes~\eqref{eq: parametric surface} to represent surfaces with a non-separable parameterization as
\begin{equation}
\boldsymbol\sigma(u,v)=\sum_{k\in\mathbb{Z}}\sum_{l\in\mathbb{Z}}\boldsymbol\sigma[k,l]\phi_{\boldsymbol\alpha_1}(u-k)\phi_{\boldsymbol\alpha_2}(v-l).
\end{equation}

{We use different families of interpolators to perfectly reproduce curves and surfaces with known parameterizations. In Section \ref{subsubsec:roman}, we detail the construction of the Roman surface. Additional examples are provided in the appendices such as the reproduction of ellipses (\ref{sec: repro ellipse}) and of the hyperbolic paraboloid (\ref{app: hyperbolic paraboloid}). The four surfaces in Figure \ref{fig:illustrations} were obtained from their classical parameterization following the same principle.}

\subsubsection{Reproduction of the Roman surface} \label{subsubsec:roman}
\noindent An illustrative example is the Roman surface whose parametrization is
\begin{align}\label{eq: roman}
\boldsymbol\sigma(u,v)&=
\begin{pmatrix}
\frac{1}{2} r^2\cos(2\pi u)\sin(4\pi v)\\
\frac{1}{2} r^2\sin(2\pi u)\sin(4\pi v)\\
r^2\cos(2\pi u)\sin(2\pi u)\cos^2(2\pi v)
\end{pmatrix}\\
&=
\begin{pmatrix}
\frac{1}{2} r^2\cos(2\pi u)\sin(4\pi v)\\
\frac{1}{2} r^2\sin(2\pi u)\sin(4\pi v)\\
\frac{1}{4}r^2\sin(4\pi u)(1+\cos(4\pi v)
\end{pmatrix}, \quad (u,v)\in \mathbb{R}^2.
\end{align}

\noindent We parameterize~\eqref{eq: roman} as a tensor-product surface of the form~\eqref{eq: parametric surface} and denote by $M_1$ and $M_2$ the number of control points related to $\phi_{\bm\alpha_1}$ and $\phi_{\bm\alpha_2}$. The surface is trigonometric in $u$ and $v$. Hence, we choose to construct the interpolators $\phi_{\boldsymbol{\alpha}_1}$ and $\phi_{\boldsymbol{\alpha}_2}$ with $\boldsymbol{\alpha}_1=\left(\frac{2\mathrm{i}\pi}{M_1},\frac{-2\mathrm{i}\pi}{M_1},\frac{4\mathrm{i}\pi}{M_1},\frac{-4\mathrm{i}\pi}{M_1}\right)$ and $\boldsymbol{\alpha}_2=\left(0,\frac{4\mathrm{i}\pi}{M_2},\frac{-4\mathrm{i}\pi}{M_2}\right)$ to express~\eqref{eq: roman} as
\begin{equation}
\boldsymbol\sigma(u,v) = \sum_{k\in \mathbb{Z}}\sum_{l\in \mathbb{Z}}\boldsymbol\sigma[k,l]\phi_{\boldsymbol{\alpha}_1}(M_1u-k)\phi_{\boldsymbol{\alpha}_2}(M_2v-l).
\end{equation}

\noindent In order to satisfy the relation $\alpha_n  - \alpha_m \neq 2 k \pi \mathrm{i}$, $k\in \mathbb{Z}$ for all pairs of distinct, purely imaginary roots, we choose $M_1=M_2=5$. To construct $\phi_{\boldsymbol{\alpha}_1}$, we see that $n_0=4$ and $N=2(n_0-1)=6$. Hence, the support of $\phi_{\boldsymbol{\alpha}_1}$ is of size $6$. Following~\eqref{eq:varphialpha}, the interpolator is expressed as
\begin{equation*}
\phi_{\boldsymbol{\alpha}_1}(t)=  \lambda[0] \beta_{\boldsymbol{\alpha}_1}(t)+\lambda[1] \bigg (\beta_{\boldsymbol{\alpha}_1}(t-\frac{1}{2})+\beta_{\boldsymbol{\alpha}_1}(t + \frac{1}{2})\bigg ) + \lambda[2] \bigg ( \beta_{\boldsymbol{\alpha}_1}(t- 1)+\beta_{\boldsymbol{\alpha}_1}(t+1) \bigg ).
\end{equation*}

\noindent By solving the corresponding system of equations~\eqref{eq:systemlambda} for the non-zero entries of $\lambda$, we find $\lambda[0]=18.118$,
 $\lambda[1]=-10.128$,
 and 
 $\lambda[2]=1.730$. For the construction of $\phi_{\boldsymbol{\alpha}_2}$, we have that $n_0=3$ and $N=2(n_0-1)=4$. The support of $\phi_{\boldsymbol{\alpha}_2}$ is therefore equal to $4$ and the interpolator is expressed as
\begin{equation*}
\phi_{\boldsymbol{\alpha}_2}(t)=  \lambda[0] \beta_{\boldsymbol{\alpha}_2}(t )+\lambda[1] \bigg (\beta_{\boldsymbol{\alpha}_2}(t-\frac{1}{2})+\beta_{\boldsymbol{\alpha}_2}(t + \frac{1}{2})\bigg ).
\end{equation*}
By solving~\eqref{eq:systemlambda}, we find that $\lambda[0]=7.396$ and $\lambda[1]=-2.825$.

\noindent Since the generator is an interpolator, the control points of the surface are given by its samples, specified by
$$
\boldsymbol\sigma(u,v)\big|_{u=k,v=l}=
\begin{pmatrix}
\frac{1}{2} r^2\cos(\frac{2\pi k}{M_1})\sin(\frac{4\pi l}{M_2})\\
\frac{1}{2} r^2\sin(\frac{2\pi k}{M_1})\sin(\frac{4\pi l}{M_2})\\
r^2\cos(\frac{2\pi k}{M_1})\sin(\frac{2\pi k}{M_1})\cos^2(\frac{2\pi l}{M_2})
\end{pmatrix}.
$$

\begin{figure}[htb]
\begin{centering}
\includegraphics[width=6cm]{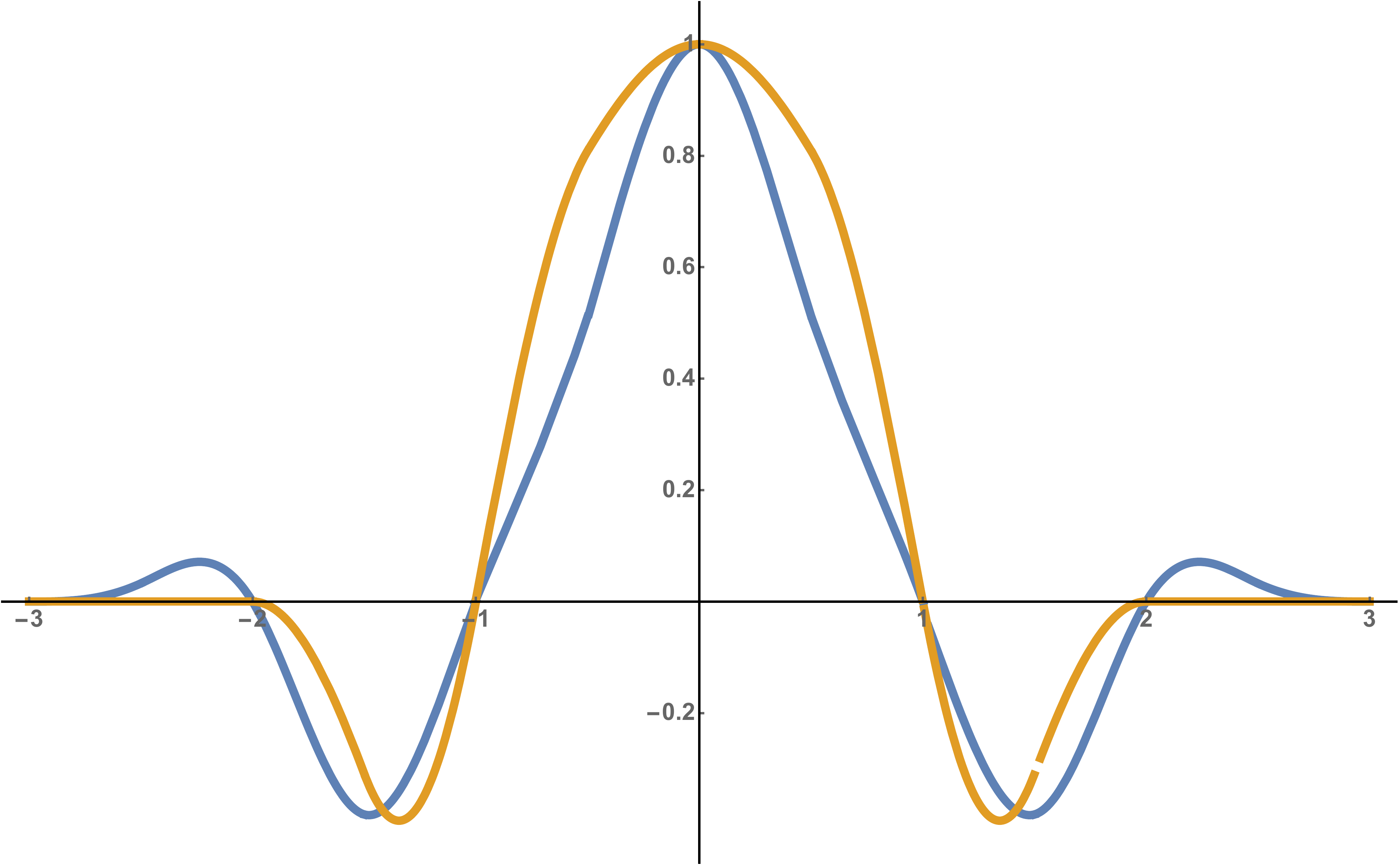}
\includegraphics[width=4cm]{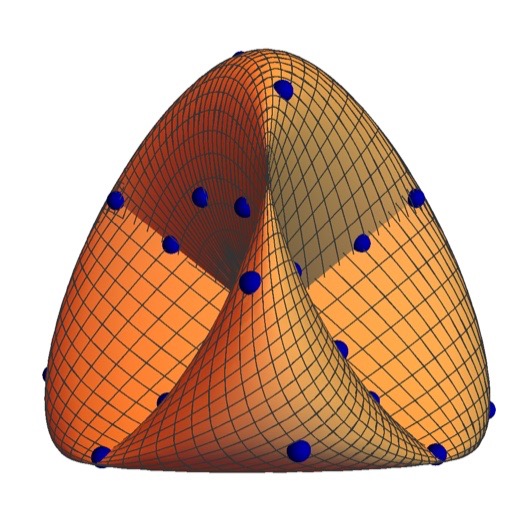}
\caption{Roman surface. The interpolators $\phi_{\boldsymbol{\alpha}_1}$ (blue) and $\phi_{\boldsymbol{\alpha}_2}$ (yellow) are shown as well as the reconstructed surface (right). The interpolatory control points are shown as blue dots on the surface.}
\label{fig: roman}
\end{centering}
\end{figure}

\noindent We choose $(u,v)\in [0,1]^2$ and $r=3$. Then, the sums in~\eqref{eq: parametric surface} are finite due to the compact support of the generators. The parameterization of the surface is given by
$
\boldsymbol\sigma(u,v)=\sum\limits_{k=-2}^{M_1+2}\sum\limits_{l=-1}^{M_2+1}\boldsymbol\sigma[k,l]\phi_{\boldsymbol{\alpha}_1}(M_1u-k)\phi_{\boldsymbol{\alpha}_2}(M_2v-l).
$
The Roman surface is illustrated in Figure~\ref{fig: roman}. 

\subsection{Interactive Shape Modeling}
\noindent The presented interpolators are well suited to be implemented in an interactive shape modeling framework; for instance, for CAD design. The key properties in such a context are 

\begin{itemize}
  \item \textit{interpolation property:} it allows to easily interact with the surface by displacing control points with a computer mouse;
  \item \textit{varying resolution:} once the ``rough'' outline of the shape is designed, the details are modeled by increasing the resolution at specific locations. 
\end{itemize}

\subsubsection{Example: Character Design}\label{sec: character design}
\noindent The interpolation property is convenient to design complex shapes as shown in Figure~\ref{fig: fancy shapes} in order to obtain a low resolution model. To increase the level of detail of the shape, we increase the resolution of the surface by first applying the pre-refinement step~\eqref{eq: modified subdiv} and then the standard refinement mask for (exponential) B-splines~\eqref{eq: refinement mask}. These two steps increase the number of control points, however, at the expense of being interpolatory. This increase in the number of control points allows one to have more flexibility in the modeling process. Furthermore, after few iterations, the convergence of the proposed modified refinement scheme allows for an interpolatory-like behavior (see Figure~\ref{fig: fancy shapes}).

\section{Discussion and Conclusion}
\noindent We have presented a general framework to construct interpolators as linear combinations of exponential B-splines of the same order $n_0$. The interpolators are compactly supported and their integer shifts form a Riesz basis whenever the corresponding B-spline does. Since the underlying building blocks are exponential B-splines, we can exploit the refinability property of the B-splines to resample the model.
Based on these general properties, we have constructed a new family of interpolators to represent parametric shapes. The new interpolators are smooth and they can be designed to perfectly reproduce polynomial, trigonometric, and hyperbolic shapes. We provide explicit examples of such generators and show in detail how idealized parametric curves and surfaces are constructed. The reconstructed shapes have the property that the control points directly lie on their boundary. This enables an intuitive manipulation of shapes by changing the location of a control point. Since the interpolators have compact support, this displacement of control points allows one to \textit{locally} control the deformation of a shape\footnote{Demo videos illustrating an implementation of our framework are found at http://bigwww.epfl.ch/demo/varying-resolution-interpolator/.}. In a next step, we plan to further investigate the refinability properties for practical applications such as real-time rendering or zooming of images.

\section{Acknowledgements}
\noindent This work was funded by the Swiss National Science Foundation under grant 200020-162343.

\appendix 
\section{Proof of Proposition~\ref{prop:riesz}} \label{prop:appendixriesz}

\begin{proof} 
We split the proof into two parts: the existence of an upper bound, relying on the one for the corresponding exponential B-spline, and the lower bound, based on the fact that the function is interpolatory. 

\paragraph{Upper Bound}
		
	We first show that one can find $B_{\bm{\alpha}}<\infty$ such that, for every $\omega \in \R$, 
	\begin{equation}\label{eq: riesz basis cond}
	\sum_{k\in \mathbb{Z}} \lvert \widehat{\beta}_{\bm{\alpha}}(\omega - 2 k \pi) \rvert^2 \leq B_{\bm{\alpha}}^2.
	\end{equation}	
	This result is well-known (see for instance~\cite[Theorem 1]{Unser0503}); we prove it for the sake of completeness.
	A more precise estimation of $B_{\bm\alpha}$ is given in~\cite[Proposition�3]{Unser0503}.
	The function $\beta_{\bm\alpha} * \beta_{\bm{\alpha}}^\vee$, where $\beta_{\bm\alpha}^\vee(t) = \beta_{\bm{\alpha}}(-t)$, is continuous and compactly supported. 
	Therefore, the sequence $c = (c[n])_{n\in \Z} = (\beta_{\bm\alpha} * \beta_{\bm{\alpha}}^\vee (n) )_{n\in \Z}$ of its samples is in $\ell_1(\Z)$. Since the Fourier transform of $\beta_{\bm\alpha} * \beta_{\bm{\alpha}}^\vee (t)$ is $\lvert \widehat{\beta}_{\bm{\alpha}} (\omega) \rvert^2$, we have that
\begin{equation}
	 \sum_{k\in \mathbb{Z}} \lvert \widehat{\beta}_{\bm{\alpha}}(\omega - 2 k \pi) \rvert^2 = \sum_{k\in \mathbb{Z}} c[k] \mathrm{e}^{- \mathrm{i} \omega k} \leq \lVert c \rVert_{\ell_1(\Z)} := B_{\bm{\alpha}}^2 < \infty.
\end{equation}
Using~\eqref{eq: interpolator fourier}, we moreover have that 
\begin{equation}
\sum_{k\in \mathbb{Z}} \lvert \widehat{\phi}_{\lambda, \bm{\alpha}}(\omega - 2 k \pi) \rvert^2
= 
\sum_{k \in \mathbb{Z}}  \left( \sum_{n\in \mathbb{Z}} \lambda[n] \mathrm{e}^{ - \mathrm{i} (\omega - 2k\pi) n /2} \right)^2 \lvert \widehat{\beta}_{\bm{\alpha}} (\omega - 2 k \pi)\rvert^2.
\end{equation} 
By splitting the sum with respect to $k$ odd or even and since $\mathrm{e}^{- \mathrm{i} (\omega- 2k\pi) n/2} = ((-1)^k)^n\mathrm{e}^{- \mathrm{i} \omega n/2}$, we have that
\begin{equation} \label{eq:equalityRieszbasis}
	 \sum_{k\in \mathbb{Z}} \lvert \widehat{\phi}_{\lambda,\bm{\alpha}}(\omega - 2 k \pi) \rvert^2  =   |G_0(\omega)|^2 \sum_{k \text{ even}}  \lvert \widehat{\beta}_{\bm{\alpha}}(\omega - 2 k \pi) \rvert^2 +   |G_1(\omega)|^2 \sum_{k \text{ odd}}  \lvert \widehat{\beta}_{\bm{\alpha}}(\omega - 2 k \pi) \rvert^2
\end{equation}
with
$	G_0(\omega) = \sum_{n \in \Z} \lambda[n] \mathrm{e}^{- \mathrm{i} \omega n/2} \text{ and }
	G_1(\omega) = \sum_{n \in \Z} (-1)^n \lambda[n] \mathrm{e}^{- \mathrm{i} \omega n/2}.
$
Clearly, for $i = 0,1$, $\lvert G_i (\omega) \rvert \leq \sum_{n\in \Z} \vert \lambda[n] \rvert = \lVert \lambda \rVert_{\ell_1(\Z)}$ and thus,
\begin{align*}
	 \sum_{k\in \mathbb{Z}} \lvert \widehat{\phi}_{\lambda,\bm{\alpha}}(\omega - 2 k \pi) \rvert^2 &\leq 
	 			\lVert \lambda \rVert_{\ell_1(\Z)}^2 \left(  \sum_{k \text{ even}}  \lvert \widehat{\beta}_{\bm{\alpha}}(\omega - 2 k \pi) \rvert^2 +    \sum_{k \text{ odd}}  \lvert \widehat{\beta}_{\bm{\alpha}}(\omega - 2 k \pi) \rvert^2 \right) \\
				& =  \lVert  \lambda \rVert_{\ell_1(\Z)}^2  \sum_{k\in \mathbb{Z}} \lvert \widehat{\beta}_{\bm{\alpha}}(\omega - 2 k \pi) \rvert^2 \\
				& \leq \lVert \lambda \rVert_{\ell_1(\Z)}^2 B_{\bm{\alpha}}^2,
\end{align*}
so that the constant $B_{\lambda, \bm\alpha} =  \lVert \lambda \rVert_{\ell_1(\Z)} B_{\bm{\alpha}}  < \infty$ acts as an upper bound in~\eqref{eq:rieszphi}. 

\paragraph{Lower Bound} 
The function $\phi_{\lambda,\bm\alpha}$ is assumed to be interpolatory; in the frequency domain, this condition is expressed as
\begin{equation} \label{eq: fourier partition of unity}
	\sum_{k \in \Z} \widehat{\phi}_{\lambda,\bm\alpha} (\omega - 2k\pi) = 1\quad\text{for all}\quad \omega \in \R.
\end{equation}
Moreover, the functions $\omega \mapsto \sum_{k\in \mathbb{Z}} \lvert \widehat{\beta}_{\bm{\alpha}}(\omega - 2 k \pi) \rvert^2$, $G_0$, and $G_1$ above are also continuous and periodic (for $G_0$ and $G_1$, this comes from $\lambda \in \ell_1(\Z)$). Therefore, the function $\omega \mapsto  \sum_{k\in \mathbb{Z}} \lvert \widehat{\phi}_{\lambda,\bm{\alpha}}(\omega - 2 k \pi) \rvert^2$ is also continuous and periodic. As such, it reaches its minimum at some frequency $\omega_0 \in [0,2\pi]$. Further, the inequality  $A_{\lambda,\bm\alpha}^2 := \sum_{k\in \mathbb{Z}} \lvert \widehat{\phi}_{\lambda,\bm{\alpha}}(\omega_0 - 2 k \pi) \rvert^2 \geq 0$ holds. Assume now that $A_{\lambda,\bm\alpha}  = 0$, then we have $\widehat{\phi}_{\bm{\alpha}}(\omega_0 - 2 k \pi) = 0$ for every $k\in \Z$, and therefore, $ \sum_{k\in \mathbb{Z}}   \widehat{\phi}_{\lambda,\bm{\alpha}}(\omega_0 - 2 k \pi)  = 0$, which contradicts~\eqref{eq: fourier partition of unity}. Hence, $A_{\lambda,\bm\alpha}  >0$ acts as a lower bound in~\eqref{eq:rieszphi}. 
\end{proof}

\noindent \textbf{Remark.} Based on~\eqref{eq:equalityRieszbasis}, we deduce the following estimates for the Riesz constants $A_{\lambda,\bm{\alpha}}$ and $B_{\lambda,\bm{\alpha}}$ associated to $\phi_{\lambda,\bm\alpha}$:
 	\begin{align}
	A_{\lambda,\bm{\alpha}} &= A_{\bm\alpha}\min_{[0,2\pi]} \lvert \hat{\lambda}(\mathrm{e}^{\mathrm i\omega})\rvert, \\
 	B_{\lambda,\bm{\alpha}} &= B_{\bm\alpha}\max_{[0,2\pi]} \lvert \hat{\lambda}(\mathrm{e}^{\mathrm i\omega})\rvert,
	\end{align}	
where $A_{\bm{\alpha}}$ and $B_{\bm\alpha}$ are the constants for the Riesz basis condition for $\beta_{\bm\alpha}$ (given in Proposition 4 and Proposition 3 in~\cite{Unser0503}), and $\hat{\lambda}(\mathrm{e}^{\mathrm i\omega}) = \sum_{n \in \Z} \lambda[n] \mathrm{e}^{- \mathrm{i} \omega n / 2}$ is the discrete Fourier transform of $\lambda$.


\section{Proof of Proposition~\ref{prop: lambda reproduction}}\label{Proof of Prop two}
\begin{proof}
The result follows from Proposition 2 in~\cite{Unser0503} which states that reproduction properties are preserved through convolution. More precisely, if $f$ is such that $\int_{-\infty}^{+\infty} f(t) \mathrm{e}^{-\alpha t} \mathrm{d} t \neq 0$ for all $\alpha \in \bm\alpha$, then $f*\beta_{\bm\alpha}$ inherits the reproduction properties of $\beta_{\bm\alpha}$. In our case, we have $\phi_{\lambda, \bm\alpha} (t) = (f*\beta_{\bm\alpha})(t)$ with $f(t) = \sum_{n\in \Z} \lambda[n] \delta(t-n/2)$. Then, for every $\alpha \in \bm\alpha$,
\begin{equation}
\int_{-\infty}^{+\infty} f(t) \mathrm{e}^{-\alpha t} \mathrm{d} t = \sum_{n\in \Z} \lambda[n] \mathrm{e}^{-\alpha n /2},
\end{equation}
which is bounded and non-zero by assumption.
\end{proof}

\section{Proof of Proposition~\ref{prop:pre filtering}}\label{Proof of Prop three}

\begin{proof}
	For the causal generator, we use~\eqref{eq:betacausalcentered} and~\eqref{eq:phi_causal} to express~\eqref{eq: interpolator fourier} as
	\begin{equation}
		\widehat{\phi}_{\lambda, \bm{\alpha}}^+ (\omega) =  \mathrm{e}^{- \mathrm{i} \omega (n_0/2 - 1)}\left( \sum_{n\in \mathbb{Z}} \lambda[n] \mathrm{e}^{- \mathrm{i} \omega n /2} \right) \widehat{\beta}_{\bm{\alpha}}^+ (\omega). 
	\end{equation}
	Then, we have
	
	\begin{align} \label{eq:phi dilation beta}
		m_0 \widehat{\phi}_{\lambda,\bm\alpha}^+ (m_0\omega) 
		&=
		\mathrm{e}^{- \mathrm{i} \omega m_0 (n_0/2 - 1)}\left( \sum_{n\in \mathbb{Z}} \lambda[n] \mathrm{e}^{- \mathrm{i} m_0 \omega n /2} \right) m_0 \widehat{\beta}_{\bm{\alpha}}^+ (m_0 \omega) 
		\nonumber \\
		&=
		 \mathrm{e}^{- \mathrm{i} \omega m_0 (n_0/2 - 1)}\left( \sum_{n\in \mathbb{Z}} \lambda[n] \mathrm{e}^{ - \mathrm{i} m_0 \omega n /2} \right) H_{\frac{\bm\alpha}{m_0},m_0}(\mathrm{e}^{\mathrm{i}\omega})  \widehat{\beta}_{\frac{\bm{\alpha}}{m_0}}^+ (\omega) \nonumber \\
		&= G_{\lambda,\frac{\bm\alpha}{m_0},m_0} (\mathrm{e}^{\mathrm{i}\omega})  \widehat{\beta}_{\frac{\bm{\alpha}}{m_0}}^+ (\omega),
	\end{align}
	where we used the relation~\eqref{eq: dilation relation} expressed in the frequency domain. Finally, we take the inverse Fourier transform of~\eqref{eq:phi dilation beta} and obtain~\eqref{eq: dilation relation for phi} in the time domain. 
\end{proof}

\section{Proof of Proposition~\ref{prop: refinement algo}}\label{Proof of Prop four}

\begin{proof}
Equation~\eqref{eq: convergence result} is equivalent to the frequency domain relation
\begin{equation}
\lim_{n\rightarrow\infty} \frac{1}{m_0 m^n} C_n(\mathrm{e}^{ \frac{\mathrm{i}\omega}{m_0m^n}}) =\widehat{f}(\omega),
\end{equation}
{where $C_n(z)=\sum_{k=-\infty}^{+\infty} c_n[k] z^{-k}$ is the $z$-transform of the discrete sequence $c_n = (c_n[k])_{k\in \mathbb{Z}}$.}
The iterative step between $c_n$ and $c_{n-1}$ in the frequency domain becomes
\begin{equation} \label{C_n C_n-1}
	C_n ( \mathrm{e}^{\frac{\mathrm{i} \omega}{m_0m^n}}) =H_{\frac{\bm\alpha}{m_0m^n}, m } ( \mathrm{e}^{\frac{\mathrm{i} \omega}{m_0m^n}}) C_{n-1}( \mathrm{e}^{\frac{\mathrm{i} \omega}{m_0 m^{n-1}}}).
\end{equation}
Iterating this relation, we obtain
\begin{equation} \label{eq : C_n C_0}
	C_n ( \mathrm{e}^{\frac{\mathrm{i} \omega}{m_0m^n}}) = \left( \prod_{k=1}^n  H_{\frac{\bm\alpha}{m_0 m^k}, m } ( \mathrm{e}^{\frac{\mathrm{i} \omega}{m_0m^k}})\right) C_0( \mathrm{e}^{\frac{\mathrm{i} \omega}{m_0}}).
\end{equation}

\noindent By expressing~\eqref{eq: dilation relation} iteratively in the frequency domain and replacing $\bm\alpha$ by $\bm\alpha / m_0$, we see that 
\begin{equation}\label{eq: mask refinement}
\begin{split}
\widehat{\beta}_{\frac{\boldsymbol\alpha}{m_0}}^+(\omega)
&=\frac{1}{m}H_{\frac{\bm\alpha}{m_0 m}, m}(\mathrm{e}^{\frac{\mathrm{i}\omega}{m}})\widehat{\beta}_{\frac{\boldsymbol\alpha}{m_0 m}}^+\left(\frac{\omega}{m} \right)  \\
&=\left( \prod_{k=1}^n\frac{1}{m}H_{\frac{\bm\alpha}{m_0 m^k}, m}(\mathrm{e}^{\frac{\mathrm{i}\omega}{m^k}}) \right) \widehat{\beta}_{\frac{\boldsymbol\alpha}{m_0 m^n}}^+\left( \frac{\omega}{m^n} \right)\\
&=\lim_{n\rightarrow\infty} \prod_{k=1}^n \frac{1}{m}H_{\frac{\bm\alpha}{m_0m^{k }}, m}(\mathrm{e}^{\frac{\mathrm{i}\omega}{m^k}}),
\end{split}
\end{equation}
where in the last line we have used the well-known convergence result from spline theory~\cite{Conti2011,Warren2001,Vonesch07}
\begin{equation} \label{eq : B spline convergence}
\lim_{n\rightarrow\infty}\widehat{\beta}_{\frac{\boldsymbol\alpha}{m_0 m^{n}}}^+ \left(\frac{\omega}{m^n} \right) = \widehat{\beta}_{\underbrace{(0,\ldots,0)}_{n_0\mathrm{times}}}(0)=\mathrm{sinc}^{n_0}(0)=1.
\end{equation}

\noindent Expressing~\eqref{eq: change of basis} in the frequency domain, we finally have
\begin{align}
	\widehat{f}(\omega) &= \frac{1}{m_0} C_0( \mathrm{e}^{\frac{\mathrm{i} \omega}{m_0}}) \widehat{\beta}_{\frac{\bm\alpha}{m_0}}^+ \left( \frac{\omega}{m_0} \right)  \nonumber \\
			&= \lim_{n \rightarrow \infty}  \frac{1}{m_0m^n}\left( \prod_{k=1}^n H_{\frac{\bm\alpha}{m_0 m^k}, m } ( \mathrm{e}^{\frac{\mathrm{i} \omega}{m_0m^k}})\right) C_0( \mathrm{e}^{\frac{\mathrm{i} \omega}{m_0}}) \\
			&= \lim_{n \rightarrow \infty} \frac{1}{m_0m^n} C_n ( \mathrm{e}^{\frac{\mathrm{i} \omega}{m_0m^n}}), \nonumber
\end{align}
where we have used~\eqref{eq : B spline convergence} and~\eqref{eq : C_n C_0} for the second and third equalities, respectively. 
\end{proof} 


\section{Reproduction of Ellipses}\label{sec: repro ellipse}
\noindent We now explicitly show how ellipses can be reproduced using our proposed interpolatory basis functions. To construct the ellipses as a function of the number of control points $M$, we choose $\boldsymbol\alpha=\left(0,\frac{2\mathrm{i}\pi}{M},-\frac{2\mathrm{i}\pi}{M}\right)$ and, hence, $n_0=3$. The interpolator is obtained by Definition~\ref{def: interpolator time domain} and by solving the corresponding system of equations~\eqref{eq:systemlambda}. The non-zero values of the sequence $\lambda$ are
\begin{equation*}
\lambda[0]=\frac{\pi ^2 \csc ^2\left(\frac{\pi }{2 M}\right) \sec \left(\frac{\pi }{M}\right)}{4
	M^2}
\end{equation*}
and
\begin{equation*}
\lambda[1]= \lambda[- 1] =-\frac{\pi ^2 \csc \left(\frac{\pi }{M}\right) \csc \left(\frac{2 \pi }{M}\right)}{M^2}.
\end{equation*}

\noindent To reproduce $\cos\left(\frac{2\pi}{M}\cdot\right)$, we take advantage of the interpolation property, which yields 
\begin{equation}
\cos\left( \frac{2\pi}{M}t \right)=
\sum\limits_{k\in\mathbb{Z}}\frac{\mathrm{e}^{\mathrm{i}\frac{2\pi}{M}k}+\mathrm{e}^{-\mathrm{i}\frac{2\pi}{M}k}}{2}\phi_{\boldsymbol\alpha}(t-k),
\end{equation}

\noindent where the coefficients are the integer samples of the curve. Normalizing the period of the cosine and using the $M$-periodized basis functions
\begin{equation}
\phi_{\boldsymbol\alpha,M}(t):=\sum\limits_{k\in \Z} \phi_{\boldsymbol\alpha}(t-Mk),
\end{equation}

\noindent we express the cosine as 
\begin{equation}
\cos(2\pi t)=\sum\limits_{k=0}^{M-1}\cos\bigg(\frac{2\pi k}{M}\bigg)\phi_{\boldsymbol\alpha,M}(Mt-k).
\end{equation}

\noindent In a similar way we obtain
\begin{equation}\label{eq: sine reproduction}
\sin(2\pi t)=\sum\limits_{k=0}^{M-1}\sin\bigg(\frac{2\pi k}{M}\bigg)\phi_{\boldsymbol\alpha,M}(Mt-k).
\end{equation}

\begin{figure}[htb]
	\begin{centering}
		\includegraphics[width=10.8cm, trim=1.5cm 0.2cm 1.5cm 3.2cm, clip=true]{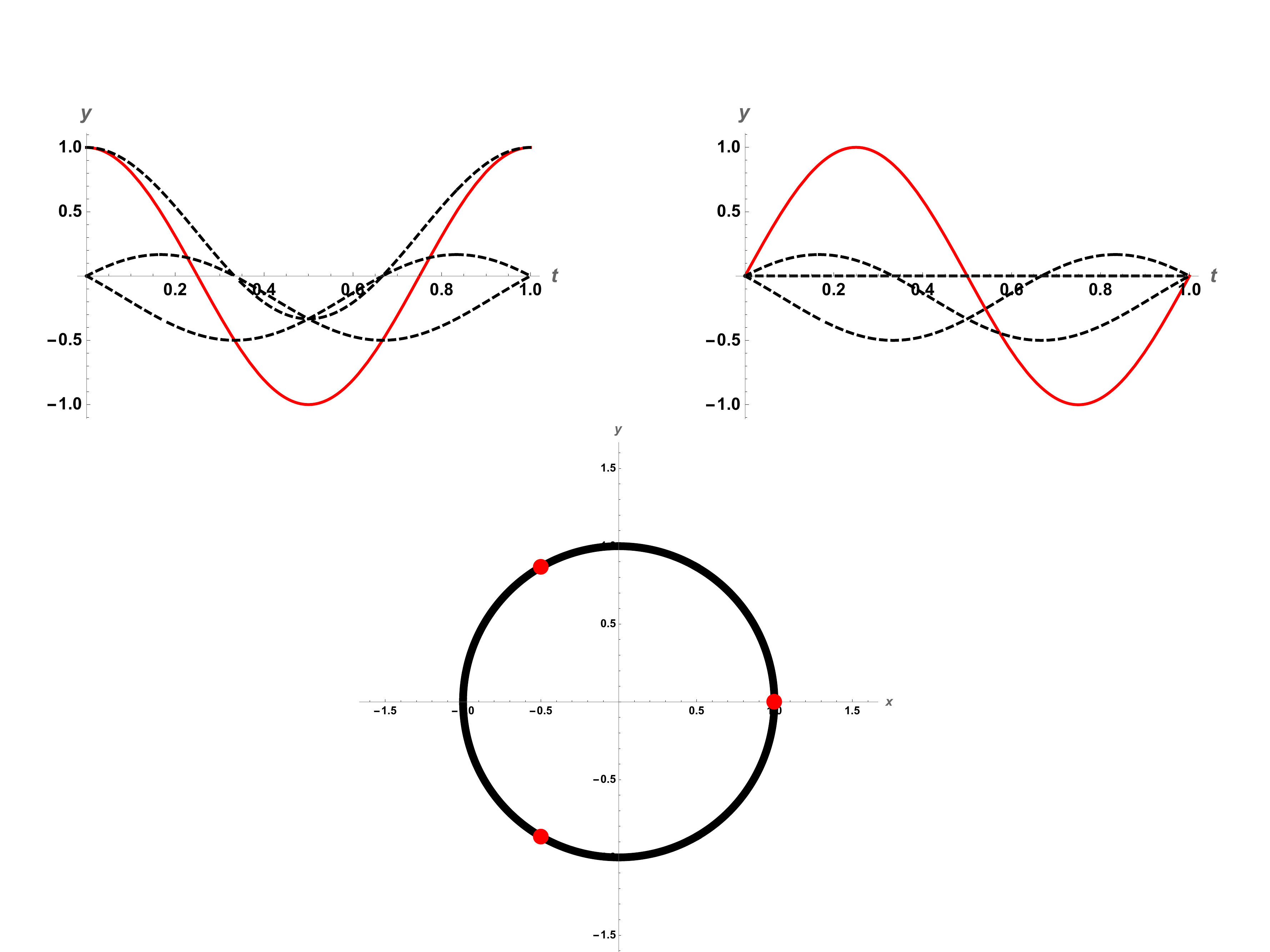}
		\caption{Top row: reproduction of the cosine (left) and sine (right) for $M=3$. The weighted and shifted basis functions are represented by dashed lines. The reconstructed parametric circle is shown in the bottom row (black) with the interpolatory control points (shown in red on the boundary of the circle). }
		\label{fig: ellipse}
	\end{centering}
\end{figure}

\noindent Plots of the trigonometric functions are shown in Figure~\ref{fig: ellipse} as well as the circle obtained through the parametric equation $\boldsymbol r(t)=(\cos(2\pi t),\sin(2\pi t))$. Ellipses can be constructed by simply applying an affine transformation to the circle $\boldsymbol r$. In order to guarantee a representation that does not depend on the location and orientation of the curve, it must be \textit{affine invariant}. This is ensured if the interpolator satisfies the partition of unity $\sum_{k\in\mathbb{Z}}\phi_{\boldsymbol\alpha,M}(\cdot-k)=1$, which implies that it must reproduce zero-degree polynomials (\textit{i.e.}, the constants). Hence, we need that $0\in\boldsymbol\alpha$.

\section{Reproduction of a Hyperbolic Paraboloid}\label{app: hyperbolic paraboloid}
\noindent A parameterization of a hyperbolic paraboloid is given by
\begin{align}\label{eq: paraboloid}
\boldsymbol\sigma(u,v)&=
\begin{pmatrix}
a u \cosh(v)\\
b u \sinh(v)\\
h u^2
\end{pmatrix}, \quad (u,v)\in \mathbb{R}^2,
\end{align}
where $a$, $b$, and $h$ are constants. The paraboloid~\eqref{eq: paraboloid} is polynomial in $u$ and hyperbolic in $v$. Hence, we choose $\boldsymbol{\alpha}_1=(0, 0,0)$ and $\boldsymbol{\alpha}_2=\left(0,\frac{1}{M_2},\frac{-1}{M_2} \right)$ when expressing~\eqref{eq: paraboloid} as the tensor-product surface
\begin{equation*}
\boldsymbol\sigma(u,v)=\sum_{k\in \mathbb{Z}}\sum_{l\in \mathbb{Z}}\boldsymbol\sigma[k,l]\phi_{\boldsymbol{\alpha}_1}(M_1u-k)\phi_{\boldsymbol{\alpha}_2}(M_2v-l).
\end{equation*}

\noindent To construct $\phi_{\boldsymbol{\alpha}_1}$, we have that $n_0=3$ and its support is equal to $N=2(n_0-1)=4$. The interpolator is expressed as
\begin{equation*}
\phi_{\boldsymbol{\alpha}_1}(t)=  \lambda[0] \beta_{\boldsymbol{\alpha}_1}(t )+\lambda[1] \bigg (\beta_{\boldsymbol{\alpha}_1}(t-\frac{1}{2})+\beta_{\boldsymbol{\alpha}_1}(t + \frac{1}{2})\bigg ).
\end{equation*}
Solving~\eqref{eq:systemlambda}, we obtain $\lambda[0]=2$ and $\lambda[1]=-\frac{1}{2}$.

\noindent For the construction of $\phi_{\boldsymbol{\alpha}_2}$, we see that $n_0=3$, $N=2(n_0-1)=4$, and its support is also of size $4$. The interpolator is given by 
\begin{equation*}
\phi_{\boldsymbol{\alpha}_2}(t)=  \lambda[0] \beta_{\boldsymbol{\alpha}_2}(t )+\lambda[1] \bigg (\beta_{\boldsymbol{\alpha}_2}(t-\frac{1}{2})+\beta_{\boldsymbol{\alpha}_2}(t + \frac{1}{2})\bigg ).
\end{equation*}
Solving~\eqref{eq:systemlambda} yields $\lambda[0]=1.968$
and 
$\lambda[1]=-0.489$. As in the previous example, the control points are obtained by sampling the surface, which leads to
\begin{equation*}
\boldsymbol\sigma(u,v)\big|_{u=k,v=l}=
\begin{pmatrix}
a \frac{k}{M_1} \cosh(\frac{l}{M_2})\\
b \frac{k}{M_1} \sinh(\frac{l}{M_2})\\
h (\frac{k}{M_1})^2
\end{pmatrix}.
\end{equation*}

\begin{figure}[htb]
	\begin{centering}
		\includegraphics[width=5.1cm,trim=0cm 0.7cm 0cm 0cm, clip=true]{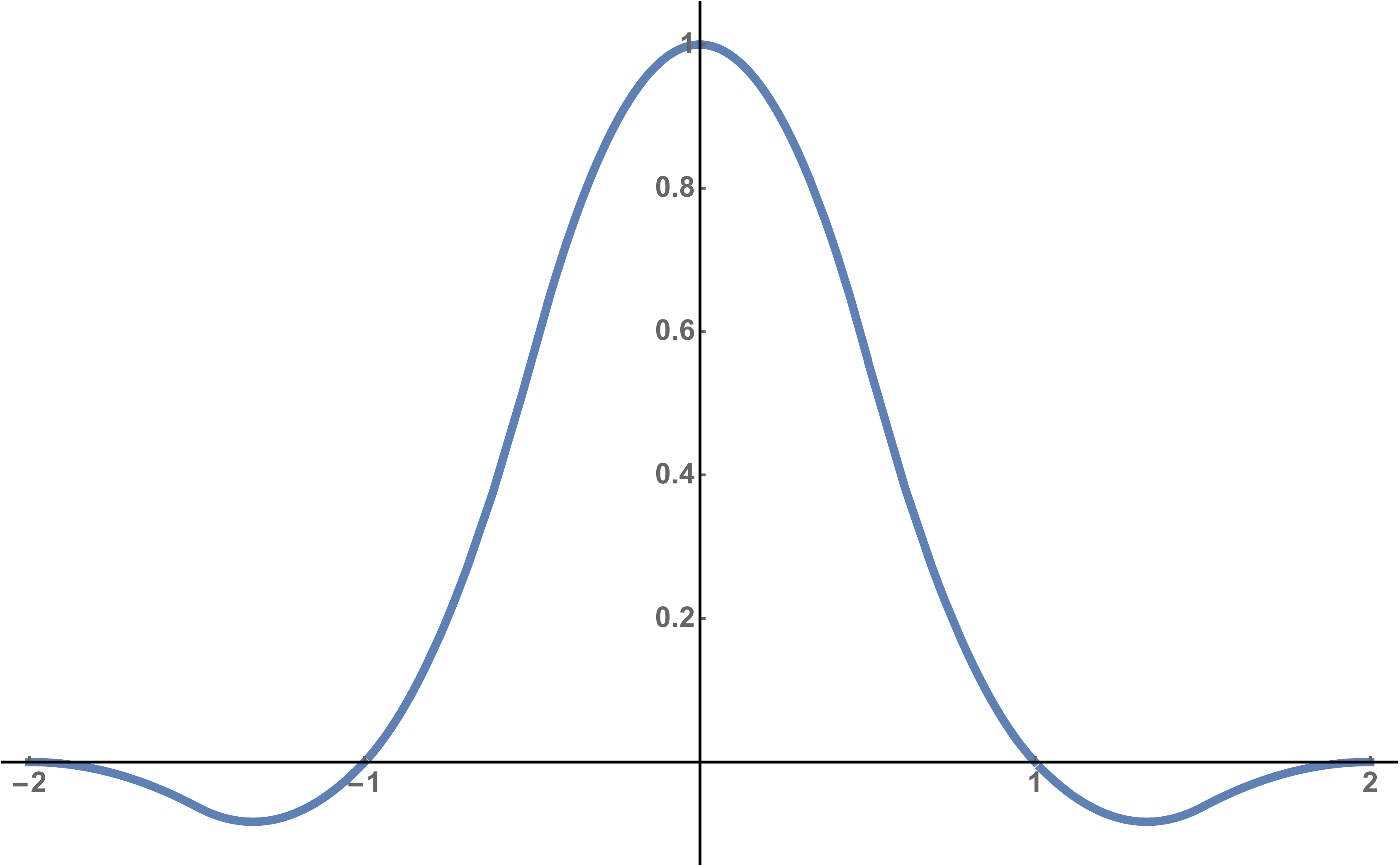}
		\includegraphics[width=5.1cm,trim=0cm 2cm 0cm 2cm, clip=true]{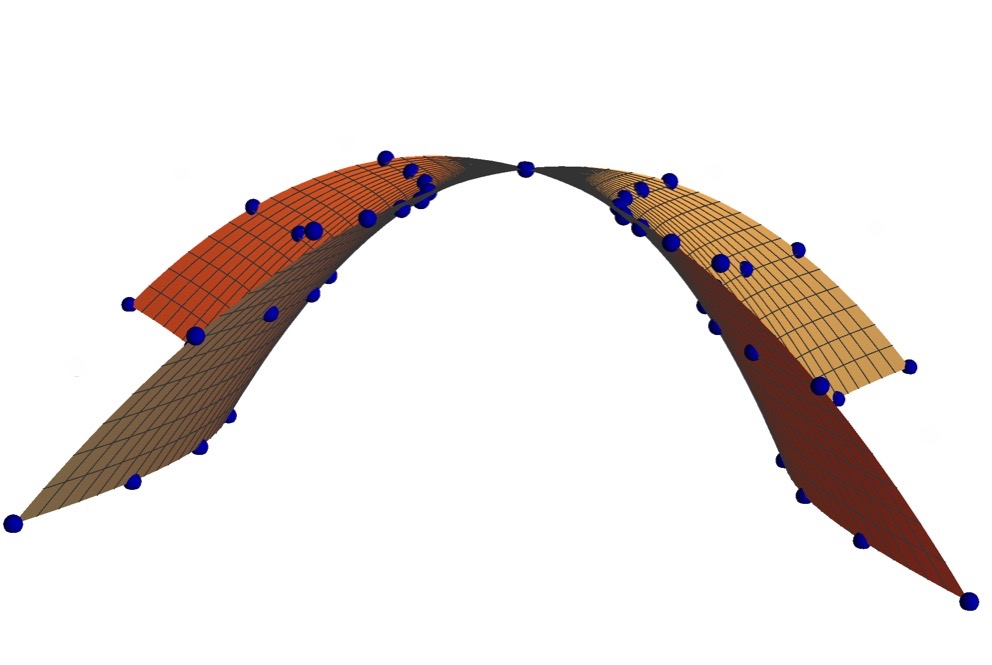}
		\caption{Hyperbolic paraboloid. On the left the interpolator $\phi_{\boldsymbol{\alpha}_2}$ is shown. ($\phi_{\boldsymbol{\alpha}_1}$ is shown in Figure~\ref{fig: interpolators}.) On the right the reconstructed hyperbolic paraboloid with its interpolatory control points (blue dots) is shown.}
		\label{fig: paraboloid}
	\end{centering}
\end{figure}

\noindent We choose $(u,v)\in [-1,1]^2$, $M_1=M_2=3$, $a=b=4$ and $h=8$. The corresponding parameterization is
\begin{equation*}
\boldsymbol\sigma(u,v)=\sum_{k=-M_1-1}^{M_1+1}\sum_{l=-M_2-1}^{M_2+1}\boldsymbol\sigma[k,l]\phi_{\boldsymbol{\alpha}_1}(M_1u-k)\phi_{\boldsymbol{\alpha}_2}(M_2v-l).
\end{equation*}
The hyperbolic paraboloid is illustrated in Figure~\ref{fig: paraboloid}. 
\section*{References}
\bibliography{Bibliography,references}

\end{document}